\documentclass[11pt]{amsart}
\usepackage[mathscr]{eucal}
\usepackage{mathrsfs}
\usepackage{amsfonts}
\usepackage{amsmath}
\usepackage{amsthm}
\usepackage{amssymb}
\usepackage{latexsym}
\usepackage[all]{xy}
\usepackage[latin1]{inputenc}
\usepackage{amscd, amssymb}
\usepackage{upgreek}
\usepackage{pstricks, pst-node}
\theoremstyle{plain}
\newtheorem{defi}{Definition}[section]
\newtheorem{thm}[defi]{Theorem}
\newtheorem{lem}[defi]{Lemma}
\newtheorem{prop}[defi]{Proposition}
\newtheorem{coro}[defi]{Corollary}
\newtheorem{rem}[defi]{Remark}

\newcommand{\pf}{\textbf{Proof.}}
\newcommand{\epf}{\hspace{14,15 cm} $\Box$}

\newcommand{\Ext}{\mathrm{Ext}}
\newcommand{\Hom}{\mathrm{Hom}}

\newcommand{\Ker}{\mathrm{Ker}}
\newcommand{\rk}{\mathrm{rk}}
\newcommand{\tr}{\mathrm{tr}}
\newcommand{\ind}{\mathrm{Ind}}
\newcommand{\dete}{\mathrm{det}}
\newcommand{\di}{\mathrm{dim}\,}

\newcommand{\Id}{\mathrm{Id}}
\newcommand{\he}{\mathrm{ht}}
\newcommand{\Su}{\mathrm{Supp}}

\newcommand{\R}{\mathbb{R}}
\newcommand{\C}{\mathbb{C}}
\newcommand{\Cc}{\hat{\mathbb{C}}}
\newcommand{\N}{\mathbb{N}}
\newcommand{\Z}{\mathbb{Z}}
\newcommand{\K}{\mathbb{K}}

\newcommand{\ga}{\mathfrak{g}}
\newcommand{\eng}{\mathcal{U}\mathfrak{g}}
\newcommand{\G}{\Gamma}
\newcommand{\GN}{\mathbf{\Gamma}_n}

\newcommand{\Pc}{\hat{\Pi}_{\lambda}(Q)}
\newcommand{\Pco}{\hat{\Pi}_{\lambda_0}}
\newcommand{\An}{\mathcal{A}_{n,\lambda,\nu}(Q)}
\newcommand{\Ann}{\mathcal{A}_{n,\lambda,\nu}}
\newcommand{\Ano}{\mathcal{A}_{n,\lambda_0, 0}}
\newcommand{\Anon}{\mathcal{A}_{n,\lambda_0+\lambda,\nu}}
\newcommand{\OGd}{{\mathcal{O}(\GN)}^{\ast}}
\newcommand{\Ogd}{{\mathcal{O}(\G^{\times n})}^{\ast}}
\newcommand{\Ogdd}{{{\mathcal{O}(\G)}^{\ast}}^{\ottimes n}}
\newcommand{\OG}{\mathcal{O}(\GN)}
\newcommand{\Og}{\mathcal{O}(\G)^{\ast}}
\newcommand{\h}{\mathcal{H}_{c, k}(\GN)}
\newcommand{\hh}{\mathcal{H}_{t, c, k}(\GN)}

\newcommand{\Ed}{\check{E}}
\newcommand{\Nn}{\mathcal{N}}
\newcommand{\Bb}{\overline{\mathcal{B}}}
\newcommand{\B}{\mathcal{B}}
\newcommand{\E}{\mathcal{E}}

\newcommand{\e}{\epsilon}
\newcommand{\ee}{\epsilon^{\ast}}
\newcommand{\jj}{\underline{j}}
\newcommand{\g}{\gamma}
\newcommand{\f}{\varphi}
\newcommand{\p}{\psi}
\newcommand{\s}{\sigma}

\newcommand{\ottimes}{\hat{\otimes}}
\newcommand{\ii}{\underline{i}}

\begin{document}

\title{Reflection functors and representations for continuous wreath-product symplectic reflection algebras}

\author{Silvia Montarani}
\address{Department of Mathematics,  Massachusetts Institute
of Technology, Cambridge, MA 02139, USA}
\email{montarani@math.mit.edu}
\maketitle
\section{\bf Introduction and main results}

Continuous symplectic reflection algebras have been recently introduced by Etingof, Gan and Ginzburg in \cite{EGG}. They are a generalization of symplectic reflection algebras (\cite{EG}) to reductive algebraic groups. 

In this paper we study representations of the wreath product continuous symplectic reflection algebra $\mathcal{H}_{c,k}(\GN)$ attached to the wreath product  $\GN:=\G^{\times n} \rtimes S_n$ of any infinite reductive subgroup $\G\subset SL(2, \C)$ with the symmetric group of rank $n$,  and to the parameter $(c, k)$,  where $k\in \C$ is a complex number, and  $c$ is an $\mathrm{Ad}(\G)$-invariant algebraic distribution on $\G$  (cfr \cite{EGG}, $\S 6$). 

In the case of a finite group $\G\subset SL(2, \C)$, a fundamental tool for the study of the wreath product symplectic reflection algebra  $\mathcal{H}_{1, c, k}(\GN)$  has been the theory of deformed preprojective algebras, introduced by Crawley-Boevey and Holland in \cite{CBH}.
 In this case, in fact, there exists a Morita equivalence  between the rank one symplectic reflection algebra $\mathcal{H}_{1,c}(\G)$ (deformed Kleinian singularity) and the deformed preprojective algebra $\Pi_{\lambda}(Q)$ attached to the (affine Dynkin) McKay quiver $Q$ of  $\G$ and to some value of the parameter $\lambda$ depending on $c$.
 
 In \cite{CBH} {\em reflection functors} were defined between the categories of finite dimensional  representations of  preprojective algebras for different values of the parameter $\lambda$. This allowed to get a complete  classification of such representations, thus of the finite dimensional representations of   $\mathcal{H}_{1,c}(\G)$. 
  Moreover, using a deformation theoretic-approach (the rank $n$ symplectic reflection algebra can be seen as a one-parameter deformation of $\Pi_{\lambda}(Q)^{\otimes n}\rtimes \C[S_n]$), it was possible to find an interesting class of finite dimensional representations for higher rank (\cite{EM}, \cite{M}).

  A second remarkable  development in the representation theory of wreath product symplectic reflection algebras of higher rank has been the introduction by Gan and Ginzburg (\cite{GG}) of the higher rank deformed preprojective algebra $\An$. The algebra $\An$ is a one-parameter deformation of the wreath product of the  preprojective algebra $\Pi_{\lambda}(Q)$ with $S_n$. In the case when $Q$ is affine Dynkin, this deformation is Morita equivalent to the higher rank symplectic reflection algebra of wreath product type. Recently,  following this interpretation of wreath product symplectic reflection algebras of higher rank in terms of deformed preprojective algebras, Gan defined a version of the reflection functors for the higher rank case (\cite{G}). This allowed him to give a more elegant and transparent formulation and proof of the  results of \cite{EM}, \cite{M}. 

In the light of these results, and using the extended McKay correspondence for reductive subgroups of $SL(2, \C)$, in this paper we define  a continuous version of  the deformed preprojective algebra $\An$, for any infinite affine Dynkin quiver $Q$  and  any rank $n$.  For appropriate values of the parameters, we prove a Morita equivalence between the continuous symplectic reflection algebra $\mathcal{H}_{c, k}(\GN)$ and the algebra $\An$, where $Q$ is the McKay quiver of $\G$. This allows us to use the results of \cite{CBH} to give a complete classification of the finite dimensional representations of the rank one continuous symplectic reflection algebra for any $\G$. For $\G=SL(2, \C)$ we compare our classification with the  results of Khare (\cite{Kh}) about the representation theory of the symplectic oscillator algebra. Moreover we show that  Gan's  definition of reflection functors easily extends to the continuous case and  we  prove the analog of the results of \cite{EM}, \cite{M} for
 finite dimensional representations in higher rank.

The structure of the paper is as follows. In Section $2$ we recall the construction of the continuous symplectic reflection algebras of wreath product type. In Section $3$ we review the theory of infinite affine Dynkin quivers. Sections $4$ and $5$ are dedicated to the definition of the algebra $\An$ and to the proof of the Morita equivalence respectively. Section $6$ contains the results about finite dimensional representations for both rank one and higher rank.

\section{{\bf Continuous symplectic reflection algebras of wreath product type}}\label{contrefalg}

\subsection{Symplectic reflections for the continuous wreath product groups $\GN$}\label{symplerefle} 
 Let $L=(\C^2, \omega_L)$,  a complex $2$-dimensional symplectic space equipped with a symplectic form $\omega_L$, and let $\G$ be an infinite reductive subgroup of $Sp(L)\cong SL(2,\C)$ (by choosing a symplectic basis). It is well known that there exist only three such groups (up to conjugation):
\begin{enumerate}
\item $SL(2,\C)$;
\item $GL(1,\C)=\C^{\ast}$, the maximal torus;
\item $\widetilde{O}_{2}$, the normalizer of the maximal torus.
\end{enumerate}
As we will see later (Section 3.1) these groups correspond to the infinite affine Dynkin diagrams $A_{+\infty}$, $A_{\infty}$, $D_{\infty}$ respectively.

We consider the space $V=L^{\oplus N}$, endowed with the induced symplectic form $\omega=\omega_L^{\oplus N}$. The symmetric group $S_n$ acts on $V$ permuting the factors. Then we have a symplectic action of the wreath product $\GN=\G^{\times n}\rtimes\C[S_n]\subset Sp(V)$ on $V$. 

For any $\g\in\G$ we will write $\g_i$ for the element $\g$ placed in the $i$th factor of $\G^{\times n}\subset\GN$. 
We will denote by  $\s_{ij}\in S_n$  the transposition $i\leftrightarrow j$.

The set $\mathcal{S}$ of symplectic reflections in $\GN$ is  defined to be  the set of elements $s$ such that  $\rk(1-s)\leq 2$ as a linear  operator on $V$. Let  $\mathcal{S}_0=\mathcal{S}\cap \G^{\times n}$ be the set of symplectic reflections in $\G^{\times n}$ and  let $\mathrm{Ad}(\GN)\s_{ij}=\left\{\s_{lm}\g_l\g_m^{-1}|l\neq m,\; \g\in\G\right\}$ be the conjugacy class of any transposition. Then it is easy to see that $\mathcal{S}=\mathcal{S}_0\cup \mathrm{Ad}(\GN)\s_{ij}$. $\GN$ acts by conjugation  on $\mathcal{S}$ preserving this decomposition. Let $\Sigma$ be the closed subscheme  of $\GN$ defined by the equation $p\circ \wedge^3(1-g|_V)=0$. $\Sigma$ is stable under conjugation. We have $\Sigma=\Sigma_0\cup \mathrm{Ad}(\GN)\s_{ij}$, where $\Sigma_0=\Sigma\cap \G^{\times n}$ and   the set of orbits $\mathcal{S}/\GN$ is a scheme isomorphic to $\Sigma/\GN$  (see \cite{EGG}, proof of Proposition $6.4$). 

\subsection{Algebraic distributions } \label{algedi}

We want here to recall some important facts about algebraic functions and distributions  for a reductive group $G$ and in particular for the groups $\G$, $\GN$. We will also introduce the notation we will use throughout the paper.

For any reductive algebraic group $G$ we will denote by $\mathcal{O}(G)$ the algebra of regular functions. 
The algebraic distributions on $G$ are the  elements of the dual space $\mathcal{O}(G)^{\ast}$. The coalgebra structure on $\mathcal{O}(G)$  induces an algebra structure on $\mathcal{O}(G)^{\ast}$ given by the convolution product. To ease notation we will write simply $\mu\mu'$ for the convolution of any two distributions $\mu, \mu'\in \mathcal{O}(G)^{\ast}$. If 
  $\Updelta:\mathcal{O}(G)\longrightarrow \mathcal{O}(G)\otimes\mathcal{O}(G)$ denotes the  the coproduct for $ \mathcal{O}(G)$, then the convolution product $\mu\mu'$ is the unique distribution on $G$ such that:
  \begin{equation}\label{convolution}
   \langle \mu\mu', f\rangle=\langle \mu\otimes\mu', \Updelta(f)\rangle \qquad \forall f\in \mathcal{O}(G)\, .
   \end{equation}

The action of $G$  on $\mathcal{O}(G)$ by left  (or right) translation induces an action on $\mathcal{O}(G)^{\ast}$ in the obvious way ($\langle\g\cdot\mu, f\rangle =\langle\mu, \g^{-1}\cdot f\rangle$, $\forall\, \g\in\G,\, \mu\in\mathcal{O}(G)^{\ast},\, f\in \mathcal{O}(G)$ ).

Similarly  the action of $G$ on itself by conjugation induces an action on $\mathcal{O}(G)$, $\mathcal{O}(G)^{\ast}$. Thus, for any closed $\mathrm{Ad}(G)$-invariant subscheme $X$ of $G$,  we have an induced action of $G$ on the space $\mathcal{O}(X)^{\ast}$ and a natural identification $\left(\mathcal{O}(X)^{\ast}\right)^{G}=\mathcal{O}(X/G)^{\ast}$. We will denote by $C(X)$ this last space.

Let us now consider our particular case. With respect to the  left (right) translation action the $\G$-modules $\mathcal{O}(\G)$, $\mathcal{O}(\G)^{\ast}$  have the following decompositions:
\begin{equation}\label{distrialgebra}
\mathcal{O}(\G)\cong\bigoplus_{i\in I}N_i\otimes N^{\ast}_{i}\cong\bigoplus_{i\in I}\mathrm{Mat}(d_i),\qquad \mathcal{O}(\G)^{\ast}\cong\prod_{i\in I}N_i\otimes N^{\ast}_{i}\cong\prod_{i\in I}\mathrm{Mat}(d_i)
\end{equation}
where $N_i$ ranges over all irreducible finite dimensional  representations of $\G$, $N_{i}^{\ast}$ denotes the dual representation, and $d_i=\dim N_i$.

In what follows we will denote by $\int_{\G}:\mathcal{O}(\G)\rightarrow \C$ the unique right and left $\G$-invariant linear form on $\mathcal{O}(\G)$  such that $\int_{\G} 1\, d\g=1$. Such form exists and is unique for any reductive group $\G$ and, if we take $N_1\cong\C$ to be the trivial representation, it can be described as the projection on the one dimensional subspace $N_1\otimes N_1^{\ast}$.

We observe now that $\mathcal{O}(S_n)=\C[S_n]$ and that the assignment  $\s\rightarrow \delta_{\s}$, where $\delta_{\s}$ denotes the delta distribution concentrated at the element $\s$, defines an algebra isomorphism  ${\mathcal{O}(S_n)}^{\ast}\cong\C[S_n]$. We have $\OG\cong \mathcal{O}(\G^{\times n})\rtimes\C[S_n]\cong \mathcal{O}(\G)^{\otimes n}\rtimes\C[S_n]$ and $\OGd~\cong~{{\mathcal{O}(\G)}^{\ast}}^{\ottimes n}\rtimes\C[S_n]$  ( where $\ottimes$ is the completed tensor product). 

We will be interested in the space $C(\mathcal{S})$ of $\GN$-invariant distributions supported on  the symplectic reflections. It follows from what we said at the end of Section \ref{symplerefle} that $C(\mathcal{S})=C(\Sigma)$ and $C(\mathcal{S})=C(\mathcal{S}_0)\oplus \C\Delta$, where $\Delta$ is the integration over $\mathrm{Ad}(\GN)\s_{ij}$ and we have a natural identification $C(\mathcal{S}_0)=C(\G)$. Thus, for any $\mu\in C(\mathcal{S})$ we can write $\mu=(c, k)$, $c\in C(\G)$, $k\in \C$.

\subsection{The wreath product construction}

In all what follows for any vector $u\in L$ we will denote by $u_i$ the corresponding vector in the $i$-th factor of $V=L^{\oplus n}$. In particular for a chosen symplectic basis $\{x, y\}$ of $L$ we will denote by $\{x_i, y_i\}_{i=1,\dots,n}$ the corresponding symplectic basis of $V$. 
For any $f\in\OG$ and $\mu\in\OGd$ we will write $f\mu=\mu f$  for the  action  of $\Og$  on $\OGd$ defined  by $\langle f\mu, g\rangle=\langle \mu, fg\rangle$.

Let $TV$ be the tensor algebra of $V$. 
\begin{defi}\label{semidire}
The semidirect product $TV\rtimes \OGd$ is the algebra generated by $w\in V$ and $\mu\in\OGd$ with the relations
$$
\mu\cdot w=\sum_i x_i\cdot (x^{\ast}_i, gw)\mu+y_i\cdot (y^{\ast}_i, gw)\mu\, ,\quad \forall w\in V,\;\mu\in\OGd
$$
where  $\left\{x^{\ast}_i, y^{\ast}_i\right\}$ denotes the dual basis of $\left\{x_i, y_i\right\}$ and $(x_i^{\ast}, gw)\mu$  denotes the action of the regular function $(x_i^{\ast}, gw)$  on $\mu$ (similarly if we substitute $x_i$ with $y_i$).
\end{defi}

\begin{rem}
Note that definition \ref{semidire} is obviously independent from the choice of the basis for $V$ as long as we choose the corresponding  dual basis for $V^{\ast}$. We use the symplectic basis $\{x_i, y_i\}$ in order not to introduce more unnecessary notation. 
\end{rem}

We will now define a family of deformations of the algebra $TV\rtimes \OGd$ indexed by  $(c,k)\in C(\mathcal{S})=C(\mathcal{S}_0)\oplus\C\Delta$.

Let  $f\in \OG\cong{{\mathcal{O}(\G)}}^{\otimes n}\rtimes\C[S_n]$ be a ``decomposable'' function, i.e.  $f=\tilde{f}(f_1\otimes\dots\otimes f_n)$, with $\tilde{f}\in \C[S_n]$ and $f_i\in\mathcal{O}(\G)$ for any $i$. Then we can write for the distribution $\Delta$:  
\begin{eqnarray}
(\Delta, f) &=&\sum_{\begin{array}{c} i,j\\i<j\end{array}}\langle\delta_{\s_{ij}},\tilde{f}\rangle\left(\int_{\G}f_i(\g)f_j({\g}^{-1}) \,d\g\,\prod_{l\neq i,j}\langle \delta_e,f_l\rangle\right)=\\ \nonumber &=& \sum_{\begin{array}{c} i,j\\i<j\end{array}}\langle\delta_{\s_{ij}},\tilde{f}\rangle\langle\Delta_{i,j}, f_1\otimes\dots\otimes f_n\rangle\nonumber
\end{eqnarray}
where  $\Delta_{ij}$ is the distribution on $\G^{\times n}$ acting as shown above and $e$ is the unit in $\G$. Thus $\Delta=\sum_{i,j|i< j}\delta_{\s_{ij}}\Delta_{ij}$.
We denote by $\omega_{L}(\g u, v)\Delta_{ij}$ the distribution on $\G^{\times n}$ such that
$$
\langle\omega_{L}(\g u, v)\Delta_{ij},  f_1\otimes\dots\otimes f_n\rangle=\int_{\G}\omega_L(\g u, v)f_i(\g)f_j({\g}^{-1}) \,d\g\,\prod_{l\neq i,j}\langle \delta_e,  f_l\rangle\, .
$$
 Finally for $c\in C(\mathcal{S})$ we will denote by $c_i$ the algebraic distribution on $\G^{\times n}$ given by $\delta_e\otimes\dots\otimes c\otimes\dots\otimes\delta_e$, where $c$ is placed in the $i$th position.

\begin{defi}\label{easyrel}
The algebra $\h$ is the quotient of $TV\rtimes\OGd$ by the following relations:
\begin{itemize}
\item[(R1)]for all $i\in[1,n]$:
\begin{equation}
[x_i, y_i]=c_i+2k\sum_{j|\,j\neq i} \delta_{\s_{ij}}\Delta_{ij}\nonumber\, .
\end{equation} 
\item[(R2)] For all $u,v\in L$, and all $i\neq j$:
\begin{equation}
[u_i,v_j]=-2k \delta_{\s_{ij}}\,\left(\omega_L(\g u, v)\Delta_{ij}\right)\nonumber
\end{equation}
\end{itemize}
\end{defi}

 In the rank one case ($n=1$) there is no parameter $k$ and, if we denote by $\{x, y\}$ a symplectic basis for the space $L$, the relations (R1), (R2)  reduce to 
\begin{equation}\label{mainre1}
[x,y]=c\, .
\end{equation}
Thus we have
\begin{equation}
\mathcal{H}_{c}(\G):=\frac{\C\langle x, y\rangle\rtimes \Og}{ \langle[x,y]-c\rangle}\, , 
\end{equation}
and the algebra $\mathcal{H}_{c}(\G)$ is the analog for continuous reductive subgroups of $SL(2,\C)$ of the deformed Kleinian singularities studied by Crawley-Boevey and Holland in \cite{CBH}. 

\begin{rem}
Note that in the original definition of Etingof, Gan and Ginzburg (\cite{EGG}, $\S$ $3.1$)  one more parameter $t\in\C$ appears. According to \cite{EGG}, if we denote  by $\delta_1$ the delta distribution at the identity element $1\in\GN$, the defining relation for the algebra $\hh$ should have the form (cfr \cite{EGG}  $\S$ $6.2$):
\begin{equation}\label{original1}
[w,w']=t\omega(w,w')\delta_{1} +\omega((1-\g)w, (1-\g)w')(c+k\Delta) \quad \forall w,w'\in V\, .
\end{equation}

In the particular case of the wreath product group $\GN=\G^{\times n}\ltimes S_n$, with $\G$  infinite, though, the parameter $t$ can be absorbed in the parameter $c$. This depends on the fact that,  if $\G\subset SL(2,\C)$ is infinite reductive, the identity element $1\in\GN$ lies in the closure of the set $\mathcal{S}$ of symplectic reflections.

Using this fact, in a similar fashion as in  the discrete case (\cite{GG}, Lemma $3.1.1$), we can reduce relation (\ref{original1}) to relations (R1), (R2) by a simple computation.

For the sake of clarity let us first look at the rank one example. In this case, according to the definition of \cite{EGG}, and absorbing the parameter $t$ in $c$, relation (\ref{mainre1}) should look like:  
\begin{equation}\label{original2}
[x,y]=\omega((1-\g)x, (1-\g)y)c
\end{equation}
from which  we  get the expression
\begin{equation}\label{original3}
[x,y]=(2-\tr_{L}(\g))c
\end{equation}
where $\tr_L$ denotes the trace in the defining representation of $\G$ on $L$.  Now it is enough to show that the invariant function $2-\tr_{L}(\g)$ is not a zero divisor in  $\mathcal{O}(\G)^{\G}$. Indeed, if this is true  then the multiplication by   $2-\tr_{L}(\g)$  is an injective linear endomorphism of  $\mathcal{O}(\G)^{\G}$, thus the induced linear endomorphism of ${\mathcal{O}(\G)^{\ast}}^{\Gamma}$ is  surjective, and any invariant distribution $c'$ can be written as $c'= (2-\tr_{L}(\g))c$, for some $c\in {\mathcal{O}(\G)^{\ast}}^{\G}=C(\G)$. But now for $\G=\C^{\ast}$ and $\G=SL(2, \C)$ the ring $O(\Gamma)^{\Gamma}$ is clearly a domain.  When $\G=\widetilde{O}_2$, the cover of the group $O_2$,  we have   $\C(\Gamma)^{\Gamma}=\C[z,z^{-1}]\oplus\C$,
where the two summands come from two connected components (so it has zero
divisors), but the function $2-\tr(\gamma)$, which is clearly not identically $0$ on the first summand,  maps to $2$ in the second summand (since
$\tr(\gamma)=0$ for $\gamma$ from the conjugacy class of orthogonal reflections),
so again it is not a zero divisor.
\footnote{Note that all what we said here is not true for continuous symplectic reflection algebras in general (see for example the case of the continuous Cherednik algebra attached to the group $O_n$, \cite {EGG} $\S$ $3.3.1$) and  the extra parameter $t$ becomes essential for a definition including all cases (\cite{EGG}, $\S$ $3.1$).}

To pass to the higher rank case, we observe first of all that $\omega(u_i, v_j)=0$ if $i\neq j$ and $\omega(x_i, y_j)=\updelta_{ij}$. Moreover, since the distribution  $\Delta$ is supported on the conjugacy class  $\mathrm{Ad}(\GN)\s_{ij}=\left\{\s_{lm}\g_l\g_m^{-1}|l\neq m,\; \g\in\G\right\}$ and for all $i\neq j$ the orbit of $\s_{ij}$ under the action of $\G^{\times n}\subset\GN$ is $\mathrm{Ad}(\G^{\times n})\s_{ij}=\left\{\s_{ij}\g_i\g_j^{-1}| \g\in\G\right\}$, we have:
$$ 
\omega((1-g)u_i, (1-g)v_j)\Delta=
$$
$$
=(1-\updelta_{ij})\left(-\omega(u_i,(\g^{-1}v)_i)-\omega((\g u)_j, v_j)\right) \delta_{\s_{ij}}\Delta_{ij}+
$$
$$
+\updelta_{ij}\sum_{l\neq i}\left(\omega(u_i, v_i)+\omega((\g u)_l, (\g v)_l)\right)\delta_{\s_{il}}\Delta_{il}=
$$
\begin{equation}\label{r2}
=-2(1-\updelta_{ij})\delta_{\s_{ij}}\left(\omega_L(\g u, v)\Delta_{ij}\right)+2\updelta_{ij}\omega_L(u,v)\sum_{l\neq j}\delta_{\s_{il}}\Delta_{il} 
\end{equation}
It's now trivial to deduce (R1), (R2) from (\ref{r2}) and the above observations.
\end{rem}

\subsection{Infinitesimal Hecke algebras}\label{infinitesimal} 
 The  rank $1$ algebra $\mathcal{H}_{c}(\G)$ has an interesting {\em infinitesimal} counterpart called the {\em infinitesimal Hecke algebra} (cfr. \cite{EGG} , Section $4$). In this section we recall the definition of such algebra.

 For $\G=SL(2,\C), GL(1,\C), \widetilde{O}_2$, let $\mathfrak{g}$ be the Lie algebra of $\G$. Then the enveloping algebra $\eng$ is naturally isomorphic to the subalgebra of $\Og$ of all algebraic distributions  set-theoretically supported at the identity element $e\in\G$ (cf. \cite{DG}, II, $\S$ $6$). 
 More precisely, if we identify any element $D\in\eng$ with the corresponding left invariant differential operator on $\G$, then the above mentioned isomorphism sends $D$ to the distribution $\tilde{D}$  such that for any function $f\in \mathcal{O}(\G)$:   
\begin{equation}\label{envedistri}
\langle\tilde{D},f\rangle:=(Df)(e)
\end{equation}
where by $Df$ we just mean $D$ applied to $f$ as a differential operator. 

In particular, $\mathrm{Ad}(\G)$-invariant distributions supported at the origin can be identified with elements of the center $\mathcal{Z}(\eng)$ of the enveloping algebra. 
If the distribution $c$ belongs to the subalgebra $\eng\subset\Og$  we define the infinitesimal Hecke algebra $\mathcal{H}_{c}(\mathfrak{g})$  as the quotient of $TV\ltimes\eng$ by the relation (\ref{mainre1}).  

When $\mathfrak{g}=\mathfrak{sl}_2$, representations of the algebra $\mathcal{H}_{c}(\mathfrak{g})$, called {\em deformed symplectic oscillator algebra} of rank $1$,  were studied by Khare in \cite{Kh}. We will compare his results with our results about finite dimensional representations of the algebra $\mathcal{H}_{c}(SL(2,\C))$ in Section \ref{Khare}.

\section{\bf Infinite quivers of affine type}\label{quiver}
\subsection{McKay correspondence for infinite reductive subgroups of $SL(2,\C)$}\label{McKay}
In this section we will briefly recall what can be called the McKay correspondence for infinite reductive subgroups of $SL(2,\C)$. Let $\G$ be as above. One can associate a graph to $\G$ in the following way. Let $\{N_i\}_{i\in I}$ (where $I=\Z$ for $\G=GL(1, \C)$ and $I=\Z_+$ for $\G=\widetilde{O}_{2}$, $SL(2, \C)$)  be  the collection of finite dimensional irreducible representations  of $\G$ and let's denote by $L$ the tautological representation of $\G$ on $\C^2$ (we recall that such a representation is self-dual). The set of vertices of the graph attached to $\G$ is indexed by $I$, while the number of edges connecting $i, j \in I$ is the multiplicity of $N_i$ in $L\otimes N_j$ (which is the same as the multiplicity of $N_j$ in $L\otimes N_i$, by self-duality of $L$). It is a classical result that  the graphs associated to $GL(1, \C)$, $\widetilde{O}_2$, $SL(2,\C)$ are the infinite  Dynkin diagrams  $A_{\infty}$, $D_{\infty}$, $A_{+\infty}$ respectively.

\begin{figure}[h]

\begin{pspicture}(0,0)(15,2 )

\dotnode(1,1){A} \dotnode(2,1){B} \dotnode(3,1){C}
\put(0,1){$\ldots$} \put(0.5,1){\line(1,0){3}}
\put(3.5,1){$\ldots$}
\put(1.7,0){$A_{\infty}$}

\dotnode(6,1){D} \dotnode(7,1){E} \dotnode(8,1){F} \dotnode(9,1){G}
\dotnode(7,2){H}
\put(6,1){\line(1,0){3.5}}
\put(9.5,1){$\ldots$}
\put(7,1){\line(0,1){1}}
\put(7.5,0){$D_\infty$}
\dotnode(12,1){I} \dotnode(13,1){J} \dotnode(14,1){K}
\put(12,1){\line(1,0){2.5}} \put(14.5,1){$\ldots$}
\put(13,0){$A_{+\infty}$}
\end{pspicture}

\caption{Graphs associated to $GL(1, \C)$, $\widetilde{O}_2$,
and $SL(2, \C)$.}
\end{figure}

Note that from the description above follows that  any vertex of such graphs corresponds to an irreducible representation $N_i$ of  $\G$, and the adjacent vertices and edges encode the decomposition of $L\otimes N_i$ into irreducibles (any such decomposition is multiplicity free, thus we get simply laced graphs). 
We will also denote by $\G$ any such graph.

\begin{rem} When $\G$ is a finite subgroup of $SL(2,\C)$, with the same procedure, one gets the affine (finite) Dynkin diagrams of type $A\, D\, E$ (standard McKay correspondence).
\end{rem}

We recall that the graphs $A_{\infty}$, $D_{\infty}$, $A_{+\infty}$, together with the analog graphs $B_{\infty}$, $C_{\infty}$, form the complete list of  connected Dynkin diagrams of infinite affine Cartan matrices, i.e. generalized Cartan matrices of infinite order, such that any principal minor of finite order is positive (\cite{K}, $\S$ $\ 4.10$). 
In particular we get the matrices:
$$ A_{\infty}=\left(\begin{array}{cccccc}\cdots &\ddots &\vdots  &\vdots & \vdots &\cdots \\
            \cdots   & -1 & 2 & -1 & \vdots  & \cdots\\
            \cdots  & \vdots  & -1 & 2 & -1 & \cdots \\
	    \cdots  &  \vdots &\vdots & \vdots &\ddots  &\cdots 
	    \end{array}\right)\qquad    
D_{\infty}=\left(\begin{array}{cccccc} 2 & 0 & -1 & \cdots &\cdots &\cdots\\
              0 & 2 & -1 & \cdots &\cdots&\cdots \\
            -1 & -1  & 2 & -1 &\vdots  &\cdots\\ 
	    0 & 0 & -1 & 2 & -1 & \cdots\\
	    \vdots & \vdots &\vdots &\vdots &\vdots &\cdots
	    \end{array}\right)
$$
$$
A_{+\infty}=\left(\begin{array}{cccc} 2 & -1 &\cdots & \cdots\\
             -1 & 2 &-1 &\cdots\\
            \vdots & \vdots &\vdots &\cdots
	    \end{array}\right)\, .
$$

In the sequel we will denote by  $A=A(\G)$ any such matrix and by $\ga'(A)$ the corresponding Kac-Moody algebra (\cite{K}, $\S$ $1$, $2$).

\subsection{Infinite rank affine root systems}\label{roots}

We want here to give a description of the root system attached to the graph $\G$ or, equivalently, to the matrix $A$.
 To this end, we observe that $A$ can be seen as the matrix of  the \emph{symmetrized Ringel form} attached to the corresponding graph, that we will denote by $(\ ,\ )$. Consider the  space $\C^I$, where $I=\Z$ for $A=A_{\infty}$, and $I=\Z_{+}$, for  $A=A_{+\infty},\, D_{\infty}$, of all column vectors $\{ \alpha_i\}_{i\in I}$, such that $\alpha_i=0$ for all but finitely many $i$. This space has a basis of coordinate vectors  $\e_i$, $i\in I$~, i.e. column vectors with $1$ at the place $i$ and $0$ elsewhere. In other words, if $\alpha\in\C^I$ then $\alpha=\sum_i\alpha_i\e_i$. The symmetrized Ringel form attached to $\G$ is the bilinear form on $\C^I$ defined as follows. Take the graph $\G$ and give any orientation to its edges.  Denote by $Q=Q(\G)$ the quiver so obtained and  by $a\in Q$ any of its arrows. If $a:i \rightarrow j$ we will say that $i=t(a)$ is the tail of $a$ and $j=h(a)$ is its head. Then we define:
\begin{equation}\label{Ringelsym}
(\alpha,\beta)=\langle \alpha, \beta\rangle + \langle \beta, \alpha  \rangle
\end{equation}
where
\begin{equation}\label{Ringel}
\langle \alpha, \beta\rangle=\sum_{i\in \Z \, (\Z_{+})} \alpha_i\,\beta_i-\sum_{a\in Q} \alpha_{t(a)}\,\beta_{h(a)}
\end{equation}

Note that  formula (\ref{Ringel})  makes sense since any of the graphs $\G$ is locally finite (i.e. any vertex has finite valency). Moreover, the matrix representing  $(\ ,\ )$ in the basis $\{\e_i\}_{i\in I}$ is exactly $A$.

We are now ready to define the root system for $A$. Our construction works more generally whenever $A$ is the matrix of  the symmetrized Ringel form for any locally finite quiver $Q$, in particular when the quiver is finite it coincides with the usual definition of root system for a quiver (\cite{CBH}, $\S$ $6$). Moreover our description coincides with the one given in \cite{K} ($\S$ $7.11$) for the root system of an infinite rank affine   Kac-Moody algebra $\ga'(A)$.  

We will say that $\e_i$ are the simple roots for $A$ (or for $Q$ or $\G$), and we will denote the set of simple roots by $\Pi=\Pi(A)$.  
Note that in our case $(\ ,\ )$ is nondegenerate on $\C^{I}$, since all the  principal minors of $A$ are. 

For any $i\in I$  we will now define the simple reflection $s_i:\C^I\rightarrow \C^I$ ($s_i:\Z^I\rightarrow \Z^I$) by:
$$
s_i(\alpha)=\alpha-(\alpha,\e_i)\e_i.
$$

The \emph{Weyl group}  $W$ attached to $A$ (equivalently to $\G$) is the group of linear automorphisms of $\C^I$ generated by the simple reflections $s_i$, $\forall i\in I$. 

 The real roots of $A$, and in general for a locally finite quiver $Q$, are defined to be the union of the orbits of the simple roots $\e_i$ under the action of $W$, we will denote the set of real roots by $\Delta^{re}=\Delta^{re}(A)$. So we have, by definition, $\Delta^{re}=\bigcup_{w\in W}w\Pi$. It is standard that any such root is positive or negative (i.e. is a sum of simple roots with all non-negative, respectively non-positive, integer coefficients) and that $\Delta^{re}_{-}=-\Delta^{re}_+$. 
 The imaginary roots of $A$, or in general for a locally finite quiver $Q$, are instead the elements of $\Z^I$ that are of the form $\pm w\beta$, for $\beta\in F$, where $F$ is the fundamental region:
$$
F=\{\beta\in \N^I \mbox{\ s.t.\ } \beta\neq 0, \mbox{\ support of\ } \beta \mbox{\ is connected},\mbox{\ and\ } (\beta,\e^i)\leq 0 \mbox{ \  } \forall i \in I\}.
$$
We denote such vectors by $\Delta^{im}=\Delta^{im}(A)$.
The root system for $A$ is the union of real and imaginary roots and we will denote it by $\Delta=\Delta(A)=\Delta^{re}\cup\Delta^{im}$. In our case, when $A$ is an infinite rank affine Cartan matrix, we have $\Delta(A)=\Delta^{re}(A)$, and there are no imaginary roots (\cite{K}, $\S$ $7.11$). This is because any infinite rank affine matrix $A$ (as well as its graph $\G$ and its root system $\Delta$),   can be seen as the limit of a sequence of positive finite rank Cartan matrices $A_n$, all of the same type, (with their Dynkin diagrams $\G_n$  and root systems $\Delta_n$), and for such matrices there are no imaginary roots (\cite{K}, $\S$ $5.2$,  Proposition $5.2$ c) ). 

All this is in contrast with the theory of finite affine Dynkin diagrams, for which the set of imaginary roots is infinite. For any finite affine Dynkin diagram, in fact,  we have  $\Delta^{im}=\{ m\delta |\, m\in\Z\}$, where $\delta=\{d_i\}$ is the vector with coordinates the dimensions $d_i=\di N_i$ of the irreducible representations attached to the vertices. The vector $\delta$ is also the  unique vector, with positive integer coordinates, of minimal norm among the vectors generating the one dimensional kernel of the corresponding affine Cartan matrix  or, equivalently, the radical of the corresponding symmetrized Ringel form, which is positive semi-definite in this case (\cite{K}, Theorem $5.6$, b) ). We can recover the analogy with the finite case if we observe that the matrix $A$ makes sense as a linear operator even on the space $\Cc^I$ of \emph{all} column  vectors $\{u_i\}$. This is because any of its rows has only finitely many non-zero entries. The kernel of this linear operator on $\Cc^I$ is one dimensional and is spanned by the vector $\mathbf{\delta}=\left\{d_i\right\}$, defined as above  for $\G$ continuous. We want to remark that this vector is {\em not} a root for the Kac-Moody algebra $\ga'(A)$ but it can be seen as  a root for a central extension of a completion of $\ga'(A)$ (see \cite{K}, $\S$ $ 7.12$).

\subsection{Action of the Weyl group on weights}

We want now to consider the dual space to $\C^I$. Such space is called the weight space and it is isomorphic to the space $\Cc^I$ above. We will denote by $\cdot$ the standard pairing between $\C^I$ and $\Cc^I$  and by $\{\ee_i\}$ the dual ``basis'' (spanning $\Cc^I$ topologically) of $\{ \e_i \}$ with respect to this pairing, that is to say $\ee_i\cdot\e_j =\updelta_{ij}$.  For any $u\in\Cc^I$ we will write $\lambda=\{\lambda_i\}$, where $\lambda=\sum_i\lambda_i\ee_i$ (where the sum  is now possibly infinite).

We will consider $\C^I$ as embedded in $\Cc^I$ via the map
$$
\begin{array}{ccccc}
\nu & : & \C^I & \longrightarrow & \Cc^I\\
    &   & \e_i & \longrightarrow & \sum_ja_{ij}\ee_j
\end{array}
$$
where $a_{ij}=(\e_i,\e_j)$. In the basis $\{\e_i\}$, $\{\ee_i\}$ the map $\nu$ is given by the Cartan matrix $A(\G)$. Moreover, for any  vector $\alpha\in\C^I$ and any $i\in I$, we have: 
\begin{equation}\label{identification}
 (\alpha,\e_i)= \nu(\alpha)\cdot \e_i=\nu(\alpha)_i\, .
\end{equation}

For any reflection $s_i$, $i\in I$, we can now consider its dual reflection  $r_i:\Cc^I\rightarrow \Cc^I $ which is uniquely determined by the property 
$$
r_i\lambda \cdot \alpha=\lambda\cdot s_i \alpha \qquad  \forall \lambda \in \Cc^I,\quad \alpha\in \C^I.
$$ 
In other words, we have $(r_i\lambda)_j=\lambda_j-(\e_i,\e_j)\lambda_i$ for any $j$, which is equivalent to: 
\begin{equation}\label{dualreflex}
r_i\lambda =\lambda-\lambda_i \nu(\e_i)\, .
\end{equation}
Thus we can define an action of $W$  on $\Cc^I$ by the condition
$$
\lambda\cdot(w\alpha) =w^{-1}\lambda\cdot \alpha \qquad \forall \lambda \in \Cc^I,\quad \alpha\in \C^I.
$$  

Let  now ``$\prec$'' be a total ordering on $\C$ satisfying the following properties (\cite{CBH}, $\S$~$7$): 
\begin{enumerate}
 \item If $a \prec b$, then $a+c\prec b+c$, for any $c\in\C$;
 \item On integers $\prec$ coincides with the usual order;
 \item For any $a\in\C$ there is $m\in\Z$ with $a\prec m$.
\end{enumerate}
An example of such an order is the lexicographic order with respect to the $\R$-basis $\{1, \sqrt{-1}\}$ of $\C$.

We say that a weight $\lambda$ is dominant if $0 \prec \lambda_i$ for all $i\in I$. Let $J\subset I$ be a finite set of indices. For any weight $\lambda$ let  $\lambda_{J}$ be the weight such that $(\lambda_{J})_j=\lambda_j$, if $j\in J$,  $(\lambda_{J})_j=0$ otherwise. Then we say that $\lambda$ is $J$-dominant if $\lambda_{J}$ is dominant. Let  $W_J$ be the  subgroup generated by the reflections $s_j$, for all $j\in J$. We have the following lemma.

\begin{lem}\label{dominance}
For any $\lambda\in \Cc^I$ and any finite subset $J\subset I$,  $\lambda$ is $W_J$-conjugate to a $J$-dominant weight $\lambda^+$.
\end{lem}

\pf\  Let us first suppose that $J$ is  connected and let $Q_J$ be the finite connected full subquiver of $Q$ corresponding to $J$. We observe that any such $Q_J$ is Dynkin. Let $J'=J\cup \partial J$, where $\partial J$ is the set of adjacent vertices for the subquiver $Q_J$ ( the vertices that are not in $Q_J$ but are joined to $Q_J$ by a path of length $1$). Let  $W(Q_{J'})$ be the Weyl group (of AD type) attached to the quiver $Q_{J'}$.   Let $U_{J'}\subset \Cc^I$ be the  vector space of weights $\mu$ satisfying $\mu_i=0$ for $i\notin J'$ (i.e. the span of $\ee_j$, $j\in J'$). As above, let $\lambda_{J'}$ be the weight  $(\lambda_{J'})_j=\lambda_j$ if $j\in J'$, $(\lambda_{J'})_j=0$ otherwise. Clearly $\lambda_{J'}\in U_{J'}$. Write $\lambda=\lambda_{J'}+(\lambda-\lambda_{J'})$. Then we have that $W_J$ fixes  $(\lambda-\lambda_{J'})$ and preserves $U_{J'}$. Identifying $U_{J'}$ with a finite dimensional vector space of dimension $|J'|$, the weight space for the finite Dynkin quiver $Q_{J'}$, we see that  $W_J$ acts on $U_{J'}$ as the parabolic subgroup $W(Q_J)\subset W(Q_{J'})$. The result now follows from the ordinary theory of Dynkin quivers. Indeed, write $\lambda_{J'}=\sum_{j\in J'} b_j\nu(\e_j)$ (this is clearly possible since any  principal minor of the Cartan matrix $A(\G)$ is non-degenerate) and  define the {\em height} of $\lambda$ as $\he(\lambda):=\sum_{j\in J'}b_j$. Consider now  a vector of the form $w\lambda_{J'}$,  for $w\in W(Q_J)$, of maximal height with respect to ``$\prec$'' (this exists since $W(Q_J)$ is finite). If $(w\lambda_{J'})_j\prec 0$ for some $j\in J$ then,  from formula (\ref{dualreflex}), we get
  $$
  \he(w\lambda_{J'})-\he(r_jw\lambda_{J'})=\he(w\lambda_{J'}-r_jw\lambda_{J'})=\he((w\lambda_{J'})_j\nu(\e_j))=(w\lambda_{J'})_j\prec 0
$$ 
 Thus $\he(w\lambda_{J'})$ is not maximal: a contradiction. So we must have  $(w\lambda_{J'})_J=(w\lambda)_J$  dominant, and  $\lambda^{+}:=w\lambda$  is  $J$-dominant.  We observe that it is possible to choose a $w$ of minimal length as an element of $W(Q_J)$ with the property $w\lambda$ is $J$-dominant. If $J$ is not connected we can clearly work separately on its connected components.

\epf

\section{\bf Continuous deformed preprojective algebras}\label{weights}

In this section we will extend the definition of deformed preprojective algebra introduced by Crawley-Boevey and Holland in \cite{CBH} to the ``continuous case'' of the affine infinite quivers of type $A_{+\infty}$, $A_{\infty}$, $D_{\infty}$. Following \cite{GG} we will then introduce  a family of (PBW) deformations of the wreath product of the newly defined \emph{continuous deformed  preprojective algebras} with the symmetric group $S_n$. It will  turn out that the  so defined algebras will provide, through a Morita equivalence, a  preprojective algebra ``interpretation'' of the continuous symplectic reflection algebras of wreath product type for any rank $n$.
\subsection{The rank 1 case}\label{rone}

 We start by recalling the Crawley-Boevey and Holland definition of the preprojective algebra.

 Let $Q$ be a quiver and $I$ the set of its vertices. Let $\overline{Q}$ be the double of $Q$, obtained by adding a reverse arrow $a^{\ast}:j\longrightarrow i$ for any arrow $a:i\longrightarrow j$ of $Q$. Let $B:=\bigoplus_{i\in I}\C$ and let $E$ be the vector space with basis  the set of edges of $\overline{Q}$. 
We have that $E$ is a $B$-bimodule and $E=\sum_{i,j\in I}E_{ij}$, where $E_{ij}$ is spanned by all the edges $a$ with $h(a)=i$, $t(a)=j$. We can form the path algebra of $\overline{Q}$ defined as $\C\overline{Q}:=T_BE=\bigoplus_{n\geq 0}T^n_BE$, where  $T^n_BE$ is the $n$-fold tensor product of $E$ over $B$. For any $i\in I$, let $e_i\in B$ denote the idempotent corresponding to the trivial path for the vertex $i$ and let 
\begin{equation}\label{preprofinite}
R_i:=\sum_{\{a\in Q|h(a)=i\}} a\cdot a^{\ast}-\sum_{\{a\in Q|t(a)=i\}} a^{\ast}\cdot a\, .
\end{equation}

If $\lambda\in B$ write $\lambda=\sum_{i\in I}\lambda_i e_i$, $\lambda_i\in \C$. 
\
\begin{defi}\label{preprofin}
 For each  $\lambda\in B$, the deformed preprojective algebra of $Q$ is the quotient
$$
\frac{\C\overline{Q}}{\langle R_i-\lambda_i e_i\rangle_{i\in I}}
$$
where $\langle\dots\rangle$ denotes the two-sided ideal generated by the indicated elements.
\end{defi}

Let now  $Q$ be any infinite locally finite quiver (in particular this is true if $Q$ is of type $A_{\infty}$, $D_{\infty}$, $A_{+\infty}$) with set of vertices $I$.  Let $B:=\Cc^I$ be the vector space over $\C$ topologically spanned by all the orthogonal idempotents $e_i$, $i\in I$ and $E$ be the vector space with topological basis formed by the edges of the double quiver $\overline{Q}$. 
Thus $E$ is a $B$-bimodule and $E=\prod_{i,j\in I}E_{ij}$, where $E_{ij}$ is spanned by all edges $a\in \overline{Q}$  such that  $h(a)=i$, $t(a)=j$. Note that, since $Q$, hence $\overline{Q}$, is locally finite, all the spaces $E_{ij}$ are finite dimensional. Then the \emph{completed} path algebra of $\overline{Q}$ is $\widehat{\C\, \overline{Q}}=\oplus_{k\geq 0}T_B^kE$, where $T_B^kE=E\otimes_B\cdots\otimes_B E$, with $k$ factors. We equip $\widehat{\C\, \overline{Q}}$ with the topology in which the family of open sets is the family of subspaces of finite codimension.  We observe that this  algebra is unital, with unit $\prod_{i\in I}e_i$, while the usual path algebra  $\C\overline{Q}\ $ (where we take  $B=\C^I$ and $E=\bigoplus E_{ij}$)  is not.

Note that, since the quiver $Q$ is locally finite,  for any $i\in I$ the element $R_i$ described in (\ref{preprofinite}) is a well defined element of  $\widehat{\C\, \overline{Q}}$.

\begin{defi}\label{prepro2}
The continuous deformed preprojective algebra $\Pc$ attached to the infinite affine quiver $Q$ and to the parameter $\lambda\in \Cc^I$ is the quotient:
$$
\Pc=\frac{\widehat{\C\, \overline{Q}}}{\langle\langle R_i-\lambda_ie_i\rangle\rangle}_{i\in I}
$$
where $\langle\langle\dots \rangle\rangle$ is the closed ideal generated by the indicated elements in the completed path algebra $\widehat{\C\, \overline{Q}}$.
\end{defi}

\subsection{Higher rank}\label{highran}
The definition of higher rank continuous deformed preprojective algebra given in this section is just a generalization to the continuous case of the one given by Gan and Ginzburg in \cite{GG}, ($1.2$). Namely, let us fix a positive integer $n$. Let $\B=B^{\ottimes n}$ be the $n$-fold completed tensor product of $B$ over $\C$. For any $l\in[1,n]$ we define the $\B$-bimodule:$$
\E_l=B^{\ottimes (l-1)}\ottimes E\ottimes B^{\ottimes (n-l)}\mbox{\ and\ } \E=\bigoplus_{1\leq l\leq n}\E_l.
$$
Note that $S_n$ acts naturally on $\E$, thus on $T_{\B}E$.
Now for any $l\in[1,n]$, any path $a\in\widehat{\C\, \overline{Q}}$ and any $\ii=(i_1,\dots, i_n)\in I^n$ we consider the elements
$$
|_{\ii} := e_{i_1}\otimes\cdots \otimes e_{i_n}\in \B$$
and
$$
a_l|_{\ii}:=e_{i_1}\otimes\cdots\otimes ae_{i_l}\otimes\cdots\otimes e_{i_n}\in T_{\B}\E_l\, .
$$

For an arrow $a\in\overline{Q}$, if   $i_l=t(a)$ we define
$$
a_l(\ii):=(i_1',\dots,i_n')\in I^n, \mbox{\ where\ } i_m'=\left\{\begin{array}{lr} i_m & \mbox{\ if\ } m\neq l\\
h(a) \,& \mbox{\ if\ } m=l\\ \end{array}\right.
$$
\begin{defi}\label{hrankprepro}
For any $\lambda\in \Cc^I$ and $\nu\in \C$, define the $\B$-algebra $\An$ to be the quotient of $T_{\B}\E\rtimes\C[S_n]$ by the following relations:
\begin{enumerate}
\item[(I)] For any $l\in[1,n]$ and $\ii=(i_1,\dots,i_n)$:
$$
(R_{i_l}-\lambda_{i_l})_l|_{\ii}=
\nu\sum_{\{m\neq l|\, i_m=i_l\}} \s_{ml}|_{\ii}\, ;
$$
\item[(II)] For any $l, m\in [1,n]$, $l\neq m$, $a,b\in\overline{Q}$ and $\ii=(i_1,\dots,i_n)$ with $i_l=t(a)$, $i_m=t(b)$:
$$
a_l|_{b_m(\ii)}b_m|_{\ii}-b_m|_{a_l(\ii)}a_l|_{\ii}=\left\{\begin{array}{lr}\nu\s_{lm}|_{\ii} &\mbox{\ if\ } b\in Q \mbox{\ and\ } a=b^{\ast}\\-\nu\s_{lm}|_{\ii} & \mbox{\ if\ } a\in Q \mbox{\ and\ } b=a^{\ast}\\ 0 & \mbox{\ else\ } \end{array}\right.
$$
\end{enumerate}
\end{defi}
For $n=1$ there is no parameter $\nu$ and $\mathcal{A}_{1,\lambda}(Q)=\Pc$, while for $N>1$ and $\nu=0$ we have $\mathcal{A}_{n,\lambda, 0}(Q)=\Pc^{\ottimes n}\rtimes\C[S_n]$. 

\section{\bf Morita equivalence}

In this section we will establish a Morita equivalence between the continuous symplectic reflection algebra $\h$ and the algebra $\An$, where $Q$ corresponds to $\G$ under the McKay correspondence, and for a suitable choice of the parameters $\lambda$, $\nu$. Our result can be considered a generalization to the case of continuous wreath-product symplectic reflection algebras of the result of \cite{GG},  $\S$ $3$ for ordinary symplectic reflection algebras of wreath product type. Our proof will follow very closely the one of \cite{GG}, but we report it in detail since  proving this  Morita equivalence is the key result of this paper from which all the results about the representation theory of  $\h$ will follow. Also these computations will help the reader to become familiar with the language of algebraic distributions.

\subsection{\bf}Before getting started we need to introduce some more notation. 
 For  any $i\in I$  we  denote by $E_{pq}^{N_i}(\g)\in\mathcal{O}(\G)$ the $(p,q)$-th matrix coefficient for the irreducible representation $N_i$, where $1\leq p,q\leq d_i=\di N_i$. Since $\G$ is reductive, we know these functions span the algebra $\mathcal{O}(\G)$. Moreover, if we take  matrix coefficients $E^{N_i}_{pq}(\g)$, $E^{N^{\ast}_j}_{lm}(\g)$ with respect to dual bases, the  following orthogonality relation holds:  
\begin{equation}\label{orthogonal}
\int_{\G}E^{N_i}_{pq}(\g)\,E^{N^{\ast}_j}_{lm}(\g) d\g=\frac{1}{d_i}\updelta_{ij}\updelta_{pl}\updelta_{qm}\, ,
\end{equation}
where  $\int_{\G}$ is the linear form described in Section \ref{algedi}.

Using the identifications (\ref{distrialgebra}) of Section $\ref{algedi}$, let  $\Ed^{N_i}_{pq}$  be the element of $\Og$ with $1$ in the $(p,q)$-th entry of the matrix for the $i$-th summand and zero elsewhere.  We can identify $\Ed^{N_i}_{pq}$ with the distribution $\langle\Ed^{N_i}_{pq}, E^{N_j}_{lm}\rangle=\updelta_{ij}\updelta_{pl}\updelta_{qm}$. Using (\ref{orthogonal}) we can  write $\langle\Ed^{N_i}_{pq}, E^{N_j}_{lm}\rangle=d_i\int_{\G}E^{N_j}_{lm}\,  E^{N_i^{\ast}}_{pq}\,d\g$. 

It is straightforward to compute that, if $\Updelta:\mathcal{O}(\G)\longrightarrow \mathcal{O}(\G)\otimes \mathcal{O}(\G)$ denotes the coproduct for $\mathcal{O}(\G)$, then $\Updelta(E^{N_i}_{pq})=\sum_{r=1}^{d_i}E^{N_i}_{pr}\otimes E^{N_i}_{rq}$. In all what follows, when there is no ambiguity, we will just omit the sum sign over repeated indices and write $\Updelta(E^{N_i}_{pq})=E^{N_i}_{pr}\otimes E^{N_i}_{rq}$. Using just the definition of convolution product (formula (\ref{convolution}), Section $\ref{algedi}$), it is now  easy to see that   $\Ed^{N_i}_{pq}\Ed^{N_j}_{lm}~=~\updelta_{ij}\updelta_{ql}\Ed^{N_i}_{pm}$, so that $\mathcal{O}(\G)^{\ast}\cong\prod_{i\in I}\mathrm{Mat}(d_i)$ is an algebra isomorphism.

We observe now  that $V\otimes\Ogdd$ is a $\Ogdd$-bimodule with right action only on the second factor and left action defined by
$$
\mu (w\otimes\mu ')=\sum_i\left(x_i\otimes(x^{\ast}_i, gw )+y_i\otimes(y^{\ast}_i, gw)\right)\mu\mu '
$$
for all $\mu,\mu'\in\Ogd$, $w\in V$, where $\left\{x_i, y_i\right\}$, $\left\{x^{\ast}_i, y^{\ast}_i\right\}$ are as in Definition \ref{semidire}. 
We have 
\begin{equation}\label{tensor}
TV\rtimes\OGd\cong T_{\Ogdd}\left(V\otimes\Ogdd\right)\rtimes\C[S_n].
\end{equation}

\subsection{\bf}Following \cite{CBH} ($\S$ $3$) and \cite{GG} ($\S$ 3.3) we will now define the idempotents $\f_i:=\Ed^{N_i}_{11}$ and $\f=\sum_i\f_i$ in the algebra $\Og$.

For the element element $\f^{\otimes n}\in\Ogdd$ we have:
$$ 
\f^{\otimes n}=\sum_{i_1,\dots ,i_n\in I} \f_{i_1}\otimes\cdots\otimes\f_{i_n}
$$
and
\begin{equation}
\sum_{i_1, p_1,\dots,i_n, p_n}(\Ed^{N_{i_1}}_{p_11}\otimes\cdots\otimes \Ed^{N_{i_n}}_{p_n1})\f^{\otimes n}(\Ed^{N_{i_1}}_{1p_1}\otimes\cdots\otimes \Ed^{N_{i_n}}_{1p_n})=\nonumber
\end{equation}
\begin{equation}\label{idempotent}
=\sum_{i_1, p_1,\dots,i_n, p_n}\Ed^{N_{i_1}}_{p_11}\Ed^{N_{i_1}}_{1p_1}\otimes\cdots\otimes \Ed^{N_{i_n}}_{p_n1}\Ed^{N_{i_n}}_{1p_n}=\delta_e^{\otimes n}\, .
\end{equation}
 
 \subsection{\bf}Since $\delta_e^{\otimes n}$ is the unit element in $\Ogdd$, 
equation (\ref{idempotent})  implies a Morita equivalence $\f^{\otimes n} \left( TV\rtimes\Ogdd\right)\f^{\otimes n}\sim TV\rtimes\Ogdd$. As a consequence we have  by (\ref{tensor})  that  $\f^{\otimes n} \left( TV\rtimes\Ogdd\right)\f^{\otimes n} $ is the tensor algebra of $\f^{\otimes n} \left(V\otimes\Ogdd\right)\f^{\otimes n}$ over $\f^{\otimes n} \Ogdd\f^{\otimes n}$.
But clearly we have an isomorphism
\begin{equation}\label{morita1}
\B\stackrel{\sim}{\longrightarrow} \f^{\otimes n}\Ogdd\f^{\otimes n}=\prod _{i_1,\dots,i_n}\C\cdot \f_{i_1}\otimes\cdots\otimes\f_{i_n}
\end{equation}
such that
$$
e_{i_1}\otimes\cdots\otimes e_{i_n}\longrightarrow \f_{i_1}\otimes\cdots\otimes\f_{i_n}\, .
$$
Moreover we have  bijections
$$
\f_i\Og\f_j\leftrightarrow \Hom_{\G}(N_i, N_j)\qquad \f_i\left(L\otimes\Og\right)\f_j\leftrightarrow \Hom_{\G}(N_i,L\otimes N_j).
$$
 Indeed we have that $\Og\f_j\cong N_j$ and  $\left(L\otimes\Og\right)\f_j\cong L\otimes N_j$.  The first is an irreducible finite dimensional representation and, since $\G$ is reductive, the last representation is completely reducible. Multiplying on the left by the idempotent $\f_i$ corresponds to projecting on the multiplicity space of the component of type $N_i$ of such representations.

Thus we have
\begin{equation}
\f^{\otimes n}\left(V\otimes{\Og}^{\ottimes n}\right)\f^{\otimes n}
=\prod_{i_1,\dots ,j_n}(\f_{i_1}\otimes\cdots\otimes\f_{i_n})(L^{\oplus n}\otimes\underbrace{\Og\ottimes\cdots\ottimes\Og}_n)(\f_{j_1}\otimes\cdots\f_{j_n})=\nonumber
\end{equation}
\begin{equation}\label{morita2}
=\bigoplus_{l=1}^{n}\prod_{i_1,\dots ,j_n}\Hom_{\G}(N_{i_1}, N_{j_1})\ottimes\cdots\ottimes\Hom_{\G}(N_{i_l},L\otimes N_{j_l})\ottimes\cdots \ottimes\Hom_{\G}(N_{i_n}, N_{j_n})\cong\E
\end{equation}
It follows from (\ref{morita1}) and (\ref{morita2}) that:
\begin{equation}\label{isom1}
\f^{\otimes n} \left(T_{\Ogdd}\left(V\otimes\Ogdd\right)\right)\f^{\otimes n}\cong T_{\B}\E
\end{equation}
and
\begin{equation}\label{isom2}
\f^{\otimes n}\left(TV\rtimes\Ogdd\right)\f^{\otimes n}=\f^{\otimes n} \left(T_{\Ogdd}\left(V\otimes\Ogdd\right)\rtimes\C[S_n]\right)\f^{\otimes n}\cong T_{\B}\E\rtimes\C[S_n] .
\end{equation}

Now  by (\ref{idempotent}) we have that the algebra $\h$ is Morita equivalent to $\f^{\otimes n}\h\f^{\otimes n}$. By (\ref{isom2}) we have that  $\f^{\otimes n}\h\f^{\otimes n}$ is isomorphic to some quotient of  $T_{\B}~\E~\rtimes~\C[S_n]$. We will show in the next theorem that for an appropriate choice of the parameters this quotient is exactly the one described in Definition \ref{hrankprepro}.

\subsection{\bf}We will need the following auxiliary lemma, which is the analog for infinite affine quivers  of Lemma $3.2$ of \cite{CBH}. Let $Q$ be a quiver attached to $\G$ (with any orientation). Let $\zeta$ be the linear map $\zeta:\C\rightarrow L\otimes L$ such that $1\rightarrow y\otimes x-x\otimes y$. 
\begin{lem}\label{morphisms}
For any arrow $a:i\rightarrow j$ in $Q$ there exist $\G$-module homomorphisms 
$$
\theta_a:N_i\rightarrow L\otimes N_j\quad\mbox{\ and\ }\quad \phi_a:N_j\rightarrow L\otimes N_i  
$$
such that for any vertex $i$
$$
\sum_{a\in Q, h(a)=i}(\Id_L\otimes \theta_a)\phi_a-\sum_{a\in Q, t(a)=i}(\Id_L\otimes \phi_a)\theta_a=-d_i(\zeta\otimes \Id_{N_i})
$$
as maps from $N_i$ to $L\otimes L\otimes N_i$, and such that
$$
(\omega_L\otimes \Id_{N_{t(a)}})(\Id_L\otimes \phi_a)\theta_a=-d_{h(a)}\Id_{N_{t(a)}}
$$
and
$$
(\omega_L\otimes \Id_{ N_{h(a)}})(\Id_L\otimes \theta_a)\phi_a=d_{t(a)}\Id_{N_{h(a)}}.
$$
Moreover the $\theta_a$, $\phi_a$  combine to give a basis for each of the spaces $\Hom_{\G}(N_i, L\otimes N_j)$.
\end{lem}

\pf\,  
In the case $\G$ is of type $A_{\infty}$ the same proof as in \cite{CBH}, Lemma $3.2$,  works without any change. For $\G$  of type $A_{+\infty}$, $D_{\infty}$ the proof goes as the one  in \cite{CBH} for  type $\tilde{D}_n, \tilde{E}_n$ if we observe that also in our case $Q$ is a (semi-infinite) tree, the $L\otimes N_i$ are all multiplicity free  and  the vector $\mathbf{\delta}$ is the unique vector in $\Cc^I$,  up to scalar multiples, such that $(\mathbf{\delta},\e_i)=0$ for all $i\in I$. 

\epf 
\subsection{\bf}
 Let $Q$ be as above and let $\chi_i$ be the character of the irreducible representation $N_i$. The following theorem is the analog of Theorem $3.5.2$ of \cite{GG}. 
\begin{thm}\label{mequivalence}
The algebra $\h$ is Morita equivalent to the algebra $\An$ for  $\nu=2k$ and $\lambda=\left\{\lambda_i\right\}$ where  $\lambda_i=\langle c, \chi_i\rangle$.
\end{thm}
\pf\, 
We give a proof for $n\geq 2$ since the proof for $n=1$ is similar and easier. Our proof  rephrases the proof of Theorem $3.5.2$ in \cite{GG} in the language of algebraic distributions on $\GN$. Using  equations (\ref{morita1})-(\ref{isom2})  and Lemma $\ref{morphisms}$ we can define an isomorphism 
$$
T_{\B}\E\rtimes \C[S_n]\rightarrow\f^{\otimes n}\left(TV\rtimes\OGd\right)\f^{\otimes n}
$$ 
$$
e_{i_1}\otimes\cdots\otimes e_{i_n}\cdot\s\rightarrow \f_{i_1}\otimes\cdots\otimes\f_{i_n}\cdot\s,
$$
$$
e_{i_1}\otimes\cdots\otimes a\otimes\cdots \otimes e_{i_n}\cdot\s\rightarrow \f_{i_1}\otimes\cdots\otimes \phi_a\otimes\cdots \otimes \f_{i_n}\cdot\s,
$$
$$
e_{i_1}\otimes\cdots\otimes a^{\ast}\otimes\cdots \otimes e_{i_n}\cdot\s\rightarrow \f_{i_1}\otimes\cdots\otimes \theta_a\otimes\cdots \otimes \f_{i_n}\cdot\s,
$$
for all $i_1,\dots,i_n\in I$, $a\in Q$. 

Let us denote by $J$ the subspace of $TV\rtimes\OGd$ spanned by relations (R1), (R2) of Definition $\ref{easyrel}$. Then the algebra $\h$ is the quotient of  $TV\rtimes\OGd$ by the two-sided ideal generated by $J$ and $\f^{\otimes n}\h\f^{\otimes n}$ is the quotient of $\f^{\otimes n}TV\rtimes\OGd\f^{\otimes n}$ by the ideal 
$$
\f^{\otimes n}\left(TV\rtimes\OGd\right) J \left(TV\rtimes\OGd\right)\f^{\otimes n}=
$$
$$
=\f^{\otimes n}\left(TV\rtimes\OGd\right)\f^{\otimes n}\Ogdd J\Ogdd\f^{\otimes n}\left(TV\rtimes\OGd\right)\f^{\otimes n}
$$
where the identity holds by  equation  (\ref{idempotent}). Our claim is that the image of the two sided ideal generated by $\f^{\otimes n}\Ogdd J\Ogdd\f^{\otimes n}$ under the above isomorphism is exactly the ideal of the defining relations  (I), (II) for $\An$.

Let's first consider the relations (R1). Notice that for any $\f\in\Og$  and $x,y\in L$ a symplectic basis, we have that in $TL\rtimes\Og$
\begin{equation}\label{R11}
\f(xy- yx)=(xy- yx)\f\, .
\end{equation}
In fact
$$
\f(xy- yx)=
$$
$$
=\left(x(x^{\ast}, \g x)\f+y(y^{\ast}, \g x)\f\right)y-\left(y(y^{\ast}, \g y)\f+x(x^{\ast}, \g y)\f\right)x=
$$
$$
=xx(x^{\ast}, \g y)(x^{\ast}, \g x)\f+xy(y^{\ast}, \g y)(x^{\ast}, \g x)\f +yx(x^{\ast}, \g y)(y^{\ast}, \g x)\f+yy(y^{\ast}, \g y)(y^{\ast}, \g x)\f+
$$
$$
-yx(x^{\ast}, \g x)(y^{\ast}, \g x)\f-yy(y^{\ast}, \g x)(y^{\ast}, \g y)\f -xy(y^{\ast}, \g x)(x^{\ast}, \g y)\f-xx(x^{\ast}, \g x)(x^{\ast}, \g y)\f=
$$
$$
=(xy-yx)\left((x^{\ast}, \g x)(y^{\ast}, \g y)-(x^{\ast}, \g y)(y^{\ast}, \g x)\right)\f=(xy-yx)(\dete \g)\,\f=
$$
$$
=(xy-yx)\f
$$
where the last identity holds since $\dete \g\equiv 1$ as a function on $SL(2, \C)$.
Also, since $c$ is a $\G$-invariant distribution,  for all $\f\in\Og$ we have 
\begin{equation}\label{R12}
\f c=c\f\in\Og. 
\end{equation}
Moreover, if for any $\f\in \Og$ we write $\f_i=\delta_e\otimes\cdots\otimes\f\otimes\cdots\otimes\delta_e\in \Ogdd$, where $\f$ is placed in the $i$th position,  we have that
\begin{equation}\label{R13}
\f_i\p_j\left(\delta_{\s_{ij}}\Delta_{ij}\right)=\left(\delta_{\s_{ij}}\Delta_{ij}\right)\f_i\p_j
\end{equation}
for any $\f, \p\in\Og$ and any $i,j\in[1,n]$.
To see that (\ref{R13}) holds it is enough to test the right and left hand side of the expression on a decomposable function $\tilde{f} (E^{N_{l_1}}_{p_1q_1}\otimes\cdots\otimes E^{N_{l_n}}_{p_nq_n})$. Suppose without loss of generality that $i=1$, $j=2$. For the right hand side we have
$$
\langle \left(\delta_{\s_{12}}\Delta_{12}\right)\f_1\p_2, \tilde{f} (E^{N_{l_1}}_{p_1q_1}\otimes\cdots\otimes E^{N_{l_n}}_{p_nq_n}\rangle=
$$
$$
=\langle\delta_{\s_{12}}, \tilde{f}\rangle\left(\langle\f, E^{N_{l_1}}_{r_1q_1}\rangle\langle\p, E^{N_{l_2}}_{r_2q_2}\rangle \int_{\G} E^{N_{d_{l_1}}}_{p_1 r_1}(\g) E^{N_{d_{l_2}}}_{p_2 r_2}(\g^{-1})\,d\g \right)\langle \delta_e, E^{N_{l_3}}_{p_3q_3}\rangle\cdots \langle \delta_e, E^{N_{l_n}}_{p_nq_n}\rangle=
$$
$$
= \tilde{f}(\s_{12})\left(\langle\f, E^{N_{l_1}}_{r_1q_1}\rangle \langle \p, E^{N_{l_2}}_{r_2q_2} \rangle \int_{\G} E^{N_{l_1}}_{p_1 r_1}(\g) E^{N_{l_2}^{\ast}}_{r_2 p_2}(\g) \,d\g\right) E^{N_{l_3}}_{p_3q_3}(e)\cdots  E^{N_{l_n}}_{p_nq_n}(e)=
$$
$$
=\updelta_{l_1l_2}\tilde{f}(\s_{12})\frac{1}{d_{l_1}}\langle\p, E^{N_{l_1}}_{p_1q_2}\rangle\langle\f, E^{N_{l_1}}_{p_2q_1}\rangle E^{N_{l_3}}_{p_3q_3}(e)\cdots  E^{N_{l_n}}_{p_nq_n}(e)\, ,
$$
while for the left hand side we have
$$
\langle \f_1\p_2\left(\delta_{\s_{12}}\Delta_{12}\right), \tilde{f} (E^{N_{l_1}}_{p_1q_1}\otimes\cdots\otimes E^{N_{l_n}}_{p_nq_n}\rangle=
$$
$$
=\langle \delta_{\s_{12}}\f_2\p_1\Delta_{12}, E^{N_{l_1}}_{p_1q_1}\otimes\cdots\otimes E^{N_{l_n}}_{p_nq_n}\rangle=
$$
$$
=\langle\delta_{\s_{12}}, \tilde{f}\rangle\left(\langle\p, E^{N_{l_1}}_{p_1r_1}\rangle\langle\f, E^{N_{l_2}}_{p_2r_2}\rangle \int_{\G} E^{N_{l_1}}_{r_1 q_1}(\g) E^{N_{l_2}}_{r_2 q_2}(\g^{-1})\,d\g\right)\langle \delta_e, E^{N_{l_3}}_{p_3q_3}\rangle\cdots \langle \delta_e, E^{N_{l_n}}_{p_nq_n}\rangle=
$$
$$
= \tilde{f}(\s_{12})\left(\langle\p, E^{N_{l_1}}_{p_1r_1}\rangle\langle\f, E^{N_{l_2}}_{p_2r_2}\rangle \int_{\G} E^{N_{l_1}}_{r_1 q_1}(\g) E^{N_{l_2}^{\ast}}_{q_2 r_2}(\g)\,d\g\right) E^{N_{l_3}}_{p_3q_3}(e)\cdots  E^{N_{l_n}}_{p_nq_n}(e)=
$$
$$
=\updelta_{l_1l_2}\tilde{f}(\s_{12})\frac{1}{d_{l_1}}\langle\p, E^{N_{l_1}}_{p_1q_2}\rangle \langle \f, E^{N_{l_1}}_{p_2q_1}\rangle E^{N_{l_3}}_{p_3q_3}(e)\cdots  E^{N_{l_n}}_{p_nq_n}(e)\, 
$$
where in both cases we used the fact that for any finite dimensional representation $N_i$, if we choose dual bases, we have $E^{N_i}_{pq}(\g^{-1})=E^{N_i^{\ast}}_{qp}(\g)$ for any $p,q=1,\dots, d_i$. Now using (\ref{R11}), (\ref{R12}), (\ref{R13}),  if we denote by $J_1$ the vector space spanned by relations (R1) we see that
$$
\f^{\otimes n}\Ogdd J_1\Ogdd\f^{\otimes n}=\f^{\otimes n}J_1\f^{\otimes n}\Ogdd\f^{\otimes n}.
$$
Then for any choice of $i_1,\dots, i_n\in I$ and $l\in[1, n]$ we have: 
$$
\f_{i_1}\otimes\cdots\otimes\f_{i_n}\cdot[x_l,y_l]=[x_l,y_l]\cdot \f_{i_1}\otimes\cdots\otimes\f_{i_n}=
$$
\begin{equation}\label{R1I}
=\f_{i_1}\otimes\cdots\otimes\frac{1}{d_{i_l}}\left(\sum _{a\in Q,\, h(a)=i_l}\phi_a\theta_a-\sum _{a\in Q,\, t(a)=i_l}\theta_a\phi_a\right)\otimes\cdots\otimes\f_{i_n}
\end{equation}
and
\begin{equation}\label{R1II}
\f_{i_1}\otimes\cdots\otimes\f_{i_n}c_l=\frac{\lambda_{i_l}}{d_{i_l}}\f_{i_1}\otimes\cdots\otimes\f_{i_n}.
\end{equation}
Indeed we can write  $\f_{i_1}\otimes\cdots\otimes\f_{i_n}c_l= \f_{i_1}\otimes\cdots\otimes\f_{i_l}c\otimes\dots\otimes\f_{i_n}$, and testing on a function $E^{N_j}_{pq}\in\mathcal{O}(\G)$ we have:
$$
\langle\f_{i_l}c, E^{N_j}_{pq}\rangle=
$$
$$
=\langle\f_{i_l}, E^{N_j}_{pr}\rangle \langle c, E^{N_j}_{rq}\rangle=
$$
$$
=\updelta_{i_lj}\updelta_{p1}\langle c, E^{N_j}_{1q}\rangle=\updelta_{i_lj}\updelta_{p1}\updelta_{q1}\langle c, E^{N_{i_l}}_{11}\rangle= 
$$
$$
=\frac{1}{d_{i_l}}\updelta_{i_lj}\updelta_{p1}\updelta_{q1}\langle c, \chi_{i_l}\rangle.
$$

The last identities follow from the fact that  $c$  is a $\G$-invariant distribution, thus a sum of duals of characters. More precisely  
$$
c=\sum_i\alpha_i\sum_{j=1}^{d_i}\Ed^{N_i}_{jj} \qquad \alpha_i\in\C\, ,
$$  
and one has 
$$
\langle c, E^{N_j}_{pq}\rangle=0 \qquad \mbox{\ if\ } p\neq q \, ,
$$
so that
$$
\langle c, \chi_{i_l}\rangle=
\langle\sum_i\alpha_i\sum_{j=1}^{d_i}\Ed^{N_i}_{jj}, \sum_{j=1}^{d_{i_l}}E^{N_{i_l}}_{jj}\rangle
=d_{i_l}\alpha_{i_l}=d_{i_l}\langle c, E^{N_{i_l}}_{11}\rangle.
$$

Now we claim that
\begin{equation}\label{R1III}
\left(\f_{i_1}\otimes\cdots\otimes\f_{i_n}\right)\delta_{\s_{lj}}\Delta_{lj}\left(\f_{i_1}\otimes\cdots\otimes\f_{i_n}\right)=\left\{\begin{array}{lr}\frac{1}{d_{i_l}}\left(\f_{i_1}\otimes\cdots\otimes\f_{i_n}\right)\delta_{\s_{lj}} & \mbox{\ if\ } i_j=i_l\\ 0 & \mbox{\ else\ }\end{array}\right.
\end{equation}
Indeed, supposing without loss of generality that $l=1$ $j=2$, and  testing the left hand side on a decomposable function   $\tilde{f}\left(E^{N_{l_1}}_{p_1q_1}\otimes\cdots\otimes E^{N_{l_n}}_{p_nq_n}\right)$  we get:
$$
\langle\left(\f_{i_1}\otimes\f_{i_2}\otimes\cdots\otimes\f_{i_n}\right)\delta_{\s_{12}}\Delta_{12}, \tilde{f}\left(E^{N_{l_1}}_{p_1q_1}\otimes\cdots\otimes E^{N_{l_n}}_{p_nq_n}\right)\rangle=
$$
$$
=\langle \delta_{\s_{12}}\left(\f_{i_2}\otimes\f_{i_1}\otimes\cdots\otimes\f_{i_n}\right)\Delta_{12}, \tilde{f}\left(E^{N_{l_1}}_{p_1q_1}\otimes\cdots\otimes E^{N_{l_n}}_{p_nq_n}\right)\rangle=
$$
$$
=\tilde{f}(\s_{12})\langle\f_{i_2}, E^{N_{l_1}}_{p_1r_1}\rangle\langle\f_{i_1}, E^{N_{l_2}}_{p_2r_2}\rangle \left(\int_{\G} E^{N_{l_1}}_{r_1q_1}(\g)E^{N_{l_2}^{\ast}}_{q_2r_2}(\g)\,d\g\right)\prod_{j\geq 3}\langle\f_{i_j}, E^{N_{l_j}}_{p_jq_j}\rangle=
$$
$$
=\tilde{f}(\s_{12})\updelta_{i_1 i_2}\prod_{j\geq 1}\updelta_{i_jl_j}\prod_{j\geq 1}\updelta_{p_j1}\prod_{j\geq 3}\updelta_{q_j1}\int_{\G}E^{N_{i_1}}_{1q_1}(\g)E^{N_{i_1}^{\ast}}_{q_21}(\g)\,d\g=
$$
$$
=\tilde{f}(\s_{12})\frac{\updelta_{i_1i_2}}{d_{i_1}}\prod_{j\geq 1}\updelta_{i_jl_j}\prod_{j\geq 1}\updelta_{p_j1}\prod_{j\geq 1}\updelta_{q_j1}=
$$
$$
=\frac{\delta_{i_1i_2}}{d_{i_1}}\langle\left(\f_{i_1}\otimes\f_{i_2}\otimes\cdots\otimes\f_{i_n}\right)\delta_{\s_{12}}, \tilde{f}\left(E^{N_{l_1}}_{p_1q_1}\otimes\cdots\otimes E^{N_{l_n}}_{p_nq_n}\right)\rangle.
$$

By (\ref{R1I}), (\ref{R1II}), (\ref{R1III}) we thus have that relations of type (R1) give us exactly the relations (I) in Definition \ref{hrankprepro}.

We will now find  the relations that are given by (R2). We will assume without loss of generality that $n=2$. First of all,  for any $u, v  \in L$ and any $\f,\p\in\mathcal{O}(\G)^{\ast}$, if $x, y$ is any  basis for $L$,  we can easily see that 
$$
\left(\f\otimes\p\right)\cdot[u_1,v_2]=
$$
\begin{equation}\label{R21}
=[x_1, x_2]\left((x^{\ast},  hu)\f\otimes(x^{\ast},  gv)\p\right)+[x_1, y_2]\left((x^{\ast},  hu))\f\otimes(y^{\ast},  gv)\p\right)+
\end{equation}
$$
+[y_1, x_2]\left((y^{\ast},  hu))\f\otimes(x^{\ast},  gv)\p\right)+[y_1, y_2]\left((y^{\ast},  hu))\f\otimes(y,  gv)\p\right)
$$
and similarly
$$
\left(\f\otimes\p\right)\delta_{\s_{12}}\omega_L(\g u,v)\Delta_{12}=
$$
\begin{equation}\label{R22}
=\delta_{\s_{12}}\omega_L(\g x,x)\Delta_{12}\left((x^{\ast}, hu)\f\otimes (x^{\ast}, gv)\p\right)+\delta_{\s_{12}}\omega_L(\g x,y)\Delta_{12}\left((x^{\ast}, hu)\f\otimes (y^{\ast}, gv)\p\right)+
\end{equation}
$$
+\delta_{\s_{12}}\omega_L(\g y,x)\Delta_{12}\left((y^{\ast}, hu)\f\otimes (x^{\ast}, gv)\p\right)+\delta_{\s_{12}}\omega_L(\g y,y)\Delta_{12}\left((y^{\ast}, hu)\f\otimes (y^{\ast}, gv)\p\right)\, .
$$
To prove this last identity we first have to introduce some more notation. For $\f, \p\in~{\mathcal{O}(\G)^{\ast}}$ we will write $\f_{h}$, $\p_{g}$ to indicate the variable with respect to which these distributions are considered (so $\f$ is a linear functional on functions in the variable $h$ etc\dots). Since the value of the distribution $\Delta_{12}$ on any function $f_1\otimes f_2 \in\mathcal{O}(\G)^{\otimes 2}$ can be written as  
$$
\langle \Delta_{12}, f_1\otimes f_2\rangle =\int_{{\scriptscriptstyle \begin{array}{c}\G\times\G \\ g=h\end{array}}} f_1(g) f_2(h^{-1})\, dg\, dh\, ,
$$ 
we will write $\Delta_{12, (g \, h)}$ to keep track of the variables.  Finally for  $E_{p_1q_1}^{N_{l_1}}\otimes E_{p_2q_2}^{N_{l_2}} \in\mathcal{O}(\G)^{\otimes 2}$ we will write 
 $$
 \Updelta(E_{p_1q_1}^{N_{l_1}}\otimes E_{p_2q_2}^{N_{l_2}})(g,h,g',h')=E_{p_1r_1}^{N_{l_1}}(g)\otimes E_{p_2r_2}^{N_{l_2}}(h)\otimes E_{r_1q_1}^{N_{l_1}}(g')\otimes E_{r_2q_2}^{N_{l_2}}(h')
 $$
 for the coproduct.
Let us now consider the decomposable function $f=\tilde{f}( E_{p_1q_1}^{N_{l_1}}\otimes E_{p_2q_2}^{N_{l_2}})$, where $\tilde{f}\in \C[S_n]$. We have 
$$
\langle \left(\f\otimes\p\right)\delta_{\s_{12}}\omega_L(\g u,v)\Delta_{12}\ ,\ f\rangle =
$$
$$
=\tilde{f}(\s_{12})\langle \left(\p_g\otimes\f_h\right)\left(\omega_L(g'u, v)\Delta_{12, (g'\, h')}\right),\,  \Updelta(E_{p_1q_1}^{N_{l_1}}\otimes E_{p_2q_2}^{N_{l_2}})(g,h,g',h')\rangle=
$$ 
$$
=\tilde{f}(\s_{12}) \langle\p_g\f_h\, ,\,  E_{p_1r_1}^{N_{l_1}}(g) E_{p_2r_2}^{N_{l_2}}(h)\rangle \langle \Delta_{12, (g'\,h')}\, ,\,  \omega_L(g'u, v)E_{r_1q_1}^{N_{l_1}}(g') E_{r_2q_2}^{N_{l_2}}(h')\rangle\, .
$$
Now making the change of variable $(g', h')=(g^{-1}\tilde{g}h, \g^{-1}\tilde{h}h)$ and using the fact that the integral is left and right translation invariant  we get  
\begin{equation}\label{long0}
\tilde{f}(\s_{12}) \langle\p_g\f_h\, ,\,  E_{p_1r_1}^{N_{l_1}}(g) E_{p_2r_2}^{N_{l_2}}(h)\rangle \langle \Delta_{12, (g'\,h')}\, ,\,  \omega_L(g^{-1}\tilde{g}hu, v)E_{r_1q_1}^{N_{l_1}}(g^{-1}\tilde{g}h) E_{r_2q_2}^{N_{l_2}}(h^{-1}\tilde{h}g)\rangle\, .
\end{equation}
We now observe that:
$$
\omega_L(g^{-1}\tilde{g}hu, v)=\omega_L(\tilde{g}hu, gv)=
$$
$$
\overbrace{(x^{\ast}, hu)(x^{\ast}, gv)\omega_L(\tilde{g}x, x)}^1+\overbrace{(x^{\ast}, hu)(y^{\ast}, gv)\omega_L(\tilde{g}x, y)}^2+
$$
\begin{equation}\label{long1}
+\overbrace{(y^{\ast}, hu)(x^{\ast}, gv)\omega_L(\tilde{g}y, x)}^{3}+\overbrace{(y^{\ast}, hu)(y^{\ast}, gv)\omega_L(\tilde{g}y, y)}^4
\end{equation}
and that 
$$
E_{r_1q_1}^{N_{l_1}}(g^{-1}\tilde{g}h) E_{r_2q_2}^{N_{l_2}}(h^{-1}\tilde{h}^{-1}g)=
$$
\begin{equation}\label{long2}
=E_{r_1s_1}^{N_{l_1}}(g^{-1})E_{s_1t_1}^{N_{l_1}}(\tilde{g})E_{t_1q_1}^{N_{l_1}}( h) E_{r_2s_2}^{N_{l_2}}(h^{-1}) E_{s_2t_2}^{N_{l_2}}(\tilde{h}^{-1}) E_{t_2q_2}^{N_{l_2}}(g)\, .
\end{equation}
Using (\ref{long1}) we can rewrite (\ref{long0}) as a sum of four terms.  If we use (\ref{long2}) to rewrite   the first of these terms, for example, we get  
$$
\tilde{f}(\s_{12})\langle\p_g\f_h\, ,\,  (x^{\ast}, hu)(x^{\ast}, gv)E_{p_1r_1}^{N_{l_1}}(g)E_{r_1s_1}^{N_{l_1}}(g^{-1})E_{t_1q_1}^{N_{l_1}}( h) E_{p_2, r_2}(h)E_{r_2s_2}^{N_{l_2}}(h^{-1})  E_{t_2q_2}^{N_{l_2}}(g)\rangle\cdot$$
$$
\cdot \langle \omega_L(\tilde{g}x, x)\Delta_{12, (\tilde{g}\, \tilde{h})}\, ,\,  E_{s_1t_1}^{N_{l_1}}(\tilde{g})E_{s_2t_2}^{N_{l_2}}(\tilde{h})\rangle=
$$
$$
=\tilde{f}(\s_{12})\langle\p_g\f_h\, ,\,  (x^{\ast}, hu)(x^{\ast}, gv)E_{p_1s_1}^{N_{l_1}}(e)E_{t_1q_1}^{N_{l_1}}( h) E_{p_2, s_2}(e)E_{t_2q_2}^{N_{l_2}}(g)\rangle\cdot
$$
$$ 
\cdot \langle \omega_L(\tilde{g}x, x)\Delta_{12, (\tilde{g}\, \tilde{h})}\, ,\,  E_{s_1t_1}^{N_{l_1}}(\tilde{g})E_{s_2t_2}^{N_{l_2}}(\tilde{h})\rangle=
$$
$$
=\tilde{f}(\s_{12})\langle\p_g\f_h\, ,\,  (x^{\ast}, hu)(x^{\ast}, gv)E_{t_1q_1}^{N_{l_1}}( h) E_{t_2q_2}^{N_{l_2}}(g)\rangle\cdot
$$
$$ 
\cdot \langle \omega_L(\tilde{g}x, x)\Delta_{12, (\tilde{g}\, \tilde{h})}\, ,\,  E_{p_1s_1}^{N_{l_1}}(e)E_{s_1t_1}^{N_{l_1}}(\tilde{g})E_{p_2, s_2}(e)E_{s_2t_2}^{N_{l_2}}(\tilde{h})\rangle=
$$
$$
=\tilde{f}(\s_{12})\langle\p_g\f_h\, ,\,  (x^{\ast}, hu)(x^{\ast}, gv)E_{t_1q_1}^{N_{l_1}}( h) E_{t_2q_2}^{N_{l_2}}(g)\rangle\cdot\langle \omega_L(\tilde{g}x, x)\Delta_{12, (\tilde{g}\, \tilde{h})}\, ,\,  E_{p_1t_1}^{N_{l_1}}(\tilde{g})E_{p_2, t_2}(\tilde{h})\rangle=
$$
$$
=\tilde{f}(\s_{12})\langle\left(\omega_L(\tilde{g}u, v)\Delta_{12, (\tilde{g}\, \tilde{h})}\right)\left((x^{\ast}, hu)\f_h\right)\left((x^{\ast}, gv)\p_g\right)\,,\,  \Updelta(E_{p_1q_1}^{N_{l_1}}\otimes E_{p_2q_2}^{N_{l_2}})(\tilde{g}, \tilde{h}, h, g)  \rangle=
$$
$$
=\langle\delta_{\s_{12}} \left((x^{\ast}, hu)\f\otimes(x^{\ast}, gu)\p\right)\omega_L(\g x,v)\Delta_{12}\ ,\ f\rangle
$$
where we just used the properties of the coproduct and counit (evaluation at the identity). It is of course possible to rewrite the remaining three terms in a similar way, so that we get exactly expression (\ref{R22}).

Now for any $i, j, k, l\in I$ an easy computation shows that, via  the identification   \\$TV~\rtimes~ {\mathcal{O}(\G)^{\ast}}^{\ottimes 2}\cong T_{{\mathcal{O}(\G)^{\ast}}^{\ottimes 2}}(V\otimes{\mathcal{O}(\G)^{\ast}}^{\ottimes 2})$, we have
$$
\left(\f_i\otimes\f_j\right)\left(\f\otimes\delta_e\right)[u_1, v_2]\left(\delta_e\otimes\p\right)\left(\f_k\otimes\f_l\right)=
$$
\begin{equation}\label{R2I}
=\left(\f_i\f\otimes\f_j\right)\left(u_1\otimes\left(\f_k\otimes\f_j\right)\right)\bigotimes\left(\f_k\otimes \f_j\right)\left(v_2\otimes\left(\f_k\otimes\p\f_l\right)\right)+
\end{equation}
$$
-\left(\f_i\otimes \f_j\right)\left(v_2\otimes\left(\f_i\otimes\p\f_l\right)\right)\bigotimes \left(\f_i\f\otimes\f_l\right)\left(u_1\otimes\left(\f_k\otimes\f_l\right)\right)
$$
where $\bigotimes$ denotes the product in $T_{{\mathcal{O}(\G)^{\ast}}^{\ottimes 2}}(V\otimes{\mathcal{O}(\G)^{\ast}}^{\ottimes 2})$,  and we can  see (\ref{R2I}) as an  identity between algebraic distributions on $\G^{\times 2}$  with values in $T^2V$. 
On the other hand we trivially  have  that 
$$
 \left(\f_i\otimes \f_j\right)\left(\f\otimes\delta_e\right)\left(\delta_{\s_{12}}\omega_L(\g u, v)\Delta_{12}\right)\left(\delta_e\otimes\p\right)\left(\f_k\otimes \f_l\right)=
$$
\begin{equation}\label{zerorel2}
=\delta_{\s_{12}} \left(\f_j\otimes\f_i\f\right) \omega_L(\g u, v)\Delta_{12}\left(\f_k\otimes \p\f_l\right)\, .
\end{equation}

As in \cite{GG} we observe now  that for any arrow $a\in\overline{Q}$ we can find distributions $\f_a, \p_a\in~\Ogdd$ and vectors $u_a, v_a\in L$ such that 
$$
\f_{t(a)}\f_a\left(u_a\otimes\f_{h(a)}\right)\neq 0\qquad \mbox{and}\qquad \f_{h(a)}\left(v_a\otimes\p_a\f_{t(a)}\right)\neq 0.
$$
Also in our case $Q$ has no loop vertices, thus we have that the spaces $\f_i\left(L\otimes{\mathcal{O}(\G)^{\ast}}^{\ottimes 2}\right)\f_j$ are at most one dimensional and for any $i, j\in I$ we have an identification:
$$
 \left(\f_i{\mathcal{O}(\G)^{\ast}}^{\ottimes 2}\otimes L\otimes {\mathcal{O}(\G)^{\ast}}^{\ottimes 2}\f_j\right)^{\G}\rightarrow\f_i\left(L\otimes{\mathcal{O}(\G)^{\ast}}^{\ottimes 2}\right)\f_j$$
 $$ 
 \alpha\otimes u\otimes\beta\rightarrow \alpha\left(u\otimes\beta\right)
$$
where  $\f_i{\mathcal{O}(\G)^{\ast}}^{\ottimes 2}\otimes L\otimes\f_j {\mathcal{O}(\G)^{\ast}}^{\ottimes 2}\cong N_i^{\ast}\otimes L\otimes N_j$ as $\G$-modules. Moreover again as in \cite{GG} we have a non degenerate $\G$-equivariant pairing
$$
\left(\f_i{\mathcal{O}(\G)^{\ast}}^{\ottimes 2}\otimes L\otimes {\mathcal{O}(\G)^{\ast}}^{\ottimes 2}\f_j\right)\bigotimes\left(\f_j{\mathcal{O}(\G)^{\ast}}^{\ottimes 2}\otimes L\otimes {\mathcal{O}(\G)^{\ast}}^{\ottimes 2}\f_i\right)\rightarrow \C
$$
$$
(\alpha\otimes u\otimes\beta)\bigotimes(\alpha'\otimes u\otimes\beta')\rightarrow (\alpha\beta')(\alpha'\beta)\omega_L(u, u'). 
$$
As a consequence, we can assume that, for any $a\in\overline{Q}$, we have  $\omega_L(u_a, v_a)=1$. Moreover $\f_{t(a)}\f\left(v_a\otimes\f_{h(a)}\right)=0$ if $\f_{h(a)}\left(v_a\otimes\p\f_{t(a)}\right)\neq 0$, and $\f_{h(a)}\left(u_a\otimes\p\f_{t(a)}\right)=0$ if $\f_{t(a)}\f\left(u_a\otimes\f_{h(a)}\right)\neq 0$.

Note that if $i\neq l$ or $j\neq k$ the expression (\ref{zerorel2}) is zero. To see this, let us evaluate the   distribution  $\left(\f_j\otimes\f_i\f\right) \omega_L(\g u, v)\Delta_{12}\left(\f_k\otimes \p\f_l\right)$ on a function $E_{p_1q_1}^{N_{l_1}}\otimes E_{p_2q_2}^{N_{l_2}}$. We have
$$
\langle \left(\f_j\otimes\f_i\f\right) \omega_L(\g u, v)\Delta_{12}\left(\f_k\otimes \p\f_l\right)\, ,\, E_{p_1q_1}^{N_{l_1}}\otimes E_{p_2q_2}^{N_{l_2}}\rangle=
$$
$$
= \langle \f_j\, ,\, E_{p_1r_1}^{N_{l_1}}\rangle\langle \f_j\f\, ,\, E_{p_2r_2}^{N_{l_2}}\rangle
\langle\omega_L(\g u, v)\Delta_{12}\, ,\, E_{r_1s_1}^{N_{l_1}}\otimes E_{r_2s_2}^{N_{l_2}}\rangle \langle \f_k\, ,\, E_{s_1q_1}^{N_{l_1}}\rangle\langle \p\f_l\, ,\, E_{s_2q_2}^{N_{l_2}}\rangle\,.
$$
Since last expression is zero if $j,k\neq l_1$ and $i, l\neq l_2$ the above distribution is identically $0$ if $j\neq k$ or $i\neq l$.

Thus, if $a, b\in \overline{Q}$ are two arrows such that $b\neq a^{\ast}$ or $a\neq b^{\ast}$ we get  from (\ref{R2I}), (\ref{zerorel2}) and (R2) that 
$$
(a\otimes h(b))(t(a)\otimes b)-(h(a)\otimes b)(a\otimes t(b))=0\, .
$$

 Suppose now $j=k$ and $i=l$. Consider an edge $a:i\rightarrow j$ in $\overline{Q}$ and suppose, for simplicity, $a\in Q$. We have an injection as an irreducible factor $\theta_a: N_i\hookrightarrow L\otimes N_j$. We can choose a basis $\mathbf{\xi}:=\left\{\xi_i\right\}$ of $L\otimes N_j=N_i\oplus\dots$ adapted to this decompositions into irreducibles
$$
\xi_1:=\f_i=E_{11}^{N_i}, \quad  \xi_1:=E_{21}^{N_i}, \quad \xi_3:=E_{31}^{N_i},\quad\dots,\quad\xi_{d_i}:=E_{d_i1}^{N_i},\quad\dots\, .
$$
On the other hand we can choose a basis $\mathbf{\mu}:=\left\{ \mu_i\right\}$ for $L\otimes N_j$ adapted to the tensor product
$$
\mu_1:=u_a\otimes\f_j=u_a\otimes E^{N_j}_{11}, \quad \mu_2=u_a\otimes E^{N_j}_{21},\quad\dots,\quad\mu_{2d_j}:=v_a\otimes E^{N_j}_{d_j1}.
$$
Let's now define the matrix $\tau=(\tau_{pq})$ by $\f\mu_q=\sum_p\tau_{pq}\xi_p$ and the matrix $\rho=(\rho_{pq})$ by $\p\xi_q=\sum_p\rho_{pq}\mu_q$. In other words we have $\tau= {}_{\xi}\f_{\mu}$, where  ${}_{\xi}\f_{\mu}=({}_{\xi_p}\f_{\mu_q})$ denotes the matrix representing  the linear map induced by $\f$ on $L\otimes N_j$ if  we choose the basis $\mathbf{\mu}$ for the domain and $\mathbf{\xi}$ for the image. Similarly we have $\rho=_{\mu}\p_{\xi}$.
Now, recalling that we are using the following identifications 
$$
N_j\stackrel{\phi_a}{\hookrightarrow}L\otimes N_j\stackrel{\vartheta_a}{\hookrightarrow}L\otimes L\otimes N_j\stackrel{\omega_L\otimes 1}{\rightarrow} N_j\, ,
$$
and that by Lemma $\ref{morphisms}$ this composition of morphisms equals $d_j\mathrm{Id}_{N_j}$,  we have that
\begin{equation}\label{nonzeromorph}
\f_{i}\f\left(u_a\otimes\f_j\right)=\tau_{11}\f_i \qquad \mbox{and}\qquad \f_j\left(v_a\otimes\p\f_i\right)=-\frac{\rho_{11}}{d_i}\f_j\, .
\end{equation}

We now claim that 
\begin{equation}\label{R2IV}
 \left(\f_j\otimes \f_i\f\right)\left(\omega_L(\g u_a, v_a)\Delta_{12}\right)\left(\f_j\otimes \p\f_i\right)=\frac{\tau_{11}\rho_{11}}{d_i} \f_j\otimes \f_i\, .
\end{equation}
First of all it's easy to see that 
$$
 \left(\f_j\otimes \f_i\f\right)\left(\omega_L(\g u_a, v_a)\Delta_{12}\right)\left(\f_j\otimes \p\f_i\right)=C\f_j\otimes \f_i\, ,
$$
where $C$ is some constant. 
To compute $C$ we will evaluate the left hand side of (\ref{R2IV}) on the  function $E^{N_{j}}_{11}\otimes E^{N_i}_{11}\in\mathcal{O}(\G)^{\otimes 2}$. We recall that we can see the functions  $E^{N_i}_{pq}(\g)$ as the matrix coefficients for the action of $\g$ on the direct factor $N_i\subset L\otimes N_j$ in the basis $\mathbf{\xi}$ and the functions $E^{L}_{rs}(\g)E^{N_j}_{pq}(\g)$ as the matrix coefficients for $\g$ on $L\otimes N_j$ in the basis $\mu$. We define the matrix $\alpha=\left\{\alpha_{pq}\right\}$ as the matrix of the change of basis $\mu_q=\sum_p\alpha_{pq}\xi_p$ and by $\tilde{\alpha}=(\tilde{\alpha}_{pq})$ its inverse. Accordingly to the previous notation we write ${}_{\xi}\f_{\xi}=({}_{\xi_q}\f_{\xi_p})$ (respectively ${}_{\xi}\p_{\xi}=({}_{\xi_q}\p_{\xi_p})$) for the matrix of the linear map $\f$ (respectively $\p$) where we chose the basis $\xi$ both for the domain and the image. 
$$
\langle\left(\f_j\otimes \f_i\f\right)\left(\omega_L(\g u_a, v_a)\Delta_{12}\right)\left(\f_j\otimes\p \f_i\right)\, ,\, E^{N_{j}}_{11}\otimes E^{N_i}_{11}\rangle=
$$
$$
=\sum_{r,p=1}^{d_j}\sum_{r',p'=1}^{d_i} \langle\Ed^{N_j}_{11}, E^{N_j}_{1r}\rangle \langle\f, E^{N_i}_{1r'}\rangle \left(\int_{\G}\omega_L(\g u_a, v_a) E^{N_j}_{rp}(\g) E^{N_i}_{r'p'}(\g^{-1})\, d\g\right) \langle\Ed^{N_j}_{11} ,E^{N_j}_{p1}\rangle \langle\p, E^{N_i}_{p'1}\rangle=
$$
$$
=\sum_{r', p'=1}^{d_i} \langle\f, E^{N_i}_{1r'}\rangle\left(\int_{\G}\omega_L(\g u_a, v_a) E^{N_j}_{11}(\g) E^{N_i}_{p'r'}(\g^{-1})\, d\g\right) \langle\p, E^{N_i}_{p'1}\rangle=
$$
$$
=\sum_{r', p'=1}^{d_i} \langle\f, E^{N_i}_{1r'}\rangle \left(\int_{\G} E^L_{11}(\g) E^{N_j}_{11}(\g) E^{N_i}_{p'r'}(\g^{-1})\, d\g\right) \langle\p, E^{N_i}_{p'1}\rangle=
$$

$$
=\frac{1}{d_i}\sum_{r', p'=1}^{d_i} \langle\f, E^{N_i}_{1r'}\rangle \langle\p, E^{N_i}_{p'1}\rangle \langle\Ed^{N_i}_{p'r'}, E^L_{11}(\g) E^{N_j}_{11}(\g)\rangle=
$$
$$
=\frac{1}{d_i}\sum_{r', p'=1}^{d_i} \langle\f, E^{N_i}_{1r'}\rangle \langle\p, E^{N_i}_{p'1}\rangle \langle\Ed^{N_i}_{p'r'}, \sum_{s,t=1}^{d_i}\tilde{\alpha}_{1s}\alpha_{t1} E^{N_i}_{st}\rangle=
$$
$$
=\frac{1}{d_i}\sum_{r', p'=1}^{d_i} \langle\f, E^{N_i}_{1r'}\rangle \langle\p, E^{N_i}_{p'1}\rangle \tilde{\alpha}_{1p'}\alpha_{r'1}=
$$
$$
=\frac{1}{d_i}\left(\sum_{p'=1}^{d_i}\langle\p, E^{N_i}_{p'1}\rangle\tilde{\alpha}_{1p'}\right) \left(\sum_{r'=1}^{d_i}\langle\f, E^{N_i}_{1r'}\rangle\alpha_{r'1}\right)=
$$
$$
=\frac{1}{d_i}\left(\sum_{p'=1}^{d_i}{}_{\xi_{p'}}\p_{\xi_{1}}\tilde{\alpha}_{1p'}\right) \left(\sum_{r'=1}^{d_i} {}_{\xi_1}\f_{\xi_{r'}}\alpha_{r'1}\right)=\frac{\rho_{11}\tau_{11}}{d_i}
$$
where the last identity holds since $\rho={}_{\mu}\p_{\xi}=\alpha^{-1}\, {}_{\xi}\p_{\xi}=\tilde{\alpha}\, {}_{\xi}\p_{\xi}$ and $\tau={}_{\xi}\f_{\mu}={}_{\xi}\p_{\xi}\, \alpha$.

Thus we have that $C=\frac{\rho_{11}a\tau_{11}}{d_i}$ and the identity (\ref{R2IV}) holds. So now taking $i=l$, $j=k$ $u=u_a, v=v_a, \f=\f_a, \p=\p_a$ in (\ref{R2I})  and  (\ref{R2IV}), and using (\ref{nonzeromorph}), we have that relation (R2) gives us exactly
$$
\left(a^{\ast}\otimes h(a)\right)\left(h(a)\otimes a\right)-\left(t(a) \otimes a\right)\left(a^{\ast}\otimes t(a)\right)=2k\delta_{\s_{12}}\left(h(a)\otimes t(a)\right)
$$
since $\tau_{11}, \rho_{11}\neq0$ in this case as observed above.
Also taking $u=u_a, v=v_a$ in (\ref{R2I})  and  (\ref{R2IV}) we have that,  if $\f_i\f\left(u_a\otimes\f_j\right)\neq 0$, then  $\f_j\left(u_a\otimes\p\f_j\right)=0$ and so $\rho_{d_j+1}=0$ (see (\ref{nonzeromorph})) and both sides of (R2) give zero. The same is true if we exchange the roles of $u_a$ and $v_a$. Thus the relations (R2) give exactly the relations (II) of Definition \ref{hrankprepro}.

\epf

\section{\bf Finite dimensional representations}\label{finitecase}
Theorem \ref{mequivalence} provides the connection we were looking for between the representation theory of quivers and deformed preprojective algebras and  the representation theory of the continuous symplectic algebra $\h$. In this section we will use the results of Crawley-Boevey and Holland (\cite{CBH}, \cite{CB}) to obtain a complete classification of the finite dimensional representations of the  rank 1 continuous deformed preprojective algebra $\Pc$ ($\mathcal{A}_{1,\lambda}(Q)$). Moreover, in higher rank, we will extend the results of Gan  (\cite{G}) to the continuous case, describing an interesting class of finite dimensional representations of $\An$ ($\h$).

\subsection{The rank 1 case}\label{finitecase1} The following easy result holds.
\begin{prop}\label{finrepr}
Any finite dimensional representation of the continuous deformed preprojective algebra $\Pc$ is a finite dimensional representation of some ordinary deformed preprojective algebra $\Pi_{\lambda_J}(Q_J)$, where $J\subset I$ is a finite subset of vertices, $Q_J$ is the corresponding full subquiver of $Q$, and $\lambda_J\in \C^J$ is the restriction of the parameter $\lambda$ to set of vertices $J$. Vice-versa, any finite dimensional representation of $\Pi_{\lambda_J}(Q_J)$ can be extended to a finite dimensional representation of $\Pc$.
\end{prop}
{\bf Proof.} By Definition \ref{prepro2} we have that a representation $M$ of $\Pc$ is a representation of $\Cc\, \overline{Q}$ with vector space $M_i=e_iM$ at the vertex $i$ and  linear maps  $M_a:M_{t(a)}\rightarrow M_{h(a)}$, $M_{a^{\ast}}:M_{h(a)}\rightarrow M_{t(a)}$ for each $a\in Q$ such that  the relation: 
$$
\sum_{\begin{array}{c}a\in Q\\ h(a)=i\end{array}}M_aM_{a^{\ast}}-\sum_{\begin{array}{c}a\in Q\\ t(a)=i\end{array}}M_{a^{\ast}}M_a=\lambda_i\Id_{M_i}
$$
holds for any $i$. But now, if $\di M<\infty$, we must have $\di M_i=\alpha_i<\infty$ for all $i$ and $\di M_i=0$ for all but finitely many $i$. Thus the representation $M$ is supported at a finite number of vertices and the result follows. 

Vice versa, suppose $ \Pi_{\lambda_J}(Q_J)$ admits a finite dimensional representation $M$, then we can clearly extend it to a representation of $\Pc$ by setting $M_i=0$ for $i\notin J$ and $M_a=M_{a^\ast}=0$ for $a\notin Q_J$. 

\epf

Let now $\N^I$ be the space of vectors $\alpha=\left\{\alpha_i\right\}$ such that $\alpha_i\in \N$ for all $i$ and $\alpha_i=0$ for all but finitely many $i$. For any $\alpha\in \N^I$ let $p(\alpha)$ be the function:
$$
p(\alpha):=1+\sum_{a\in Q}\alpha_{t(a)}\alpha_{h(a)}-\sum_{i\in I}\alpha_i^2.
$$
Proposition \ref{finrepr} implies the next Corollary.
\begin{coro}\label{CBHcont}
For $\lambda\in \Cc^I$ and  $\alpha\in\N^I$ the following are equivalent:
\begin{enumerate}
\item There is a simple representation of $\Pc$ of dimension vector $\alpha$
\item $\alpha$ is a positive root for some  finite, connected, full subquiver $Q_J$ of $Q$ such that $\lambda_J\alpha=0$ and $p(\alpha)\geq \sum_{k=1}^np(\beta^{(k)})$ for any decomposition $\alpha=\beta^{(1)}+\dots+\beta^{(n)}$ with $r\geq 2$ and $\beta^{(k)}$ a positive root with $\lambda\cdot\beta^{(k)}=0$ for all $k$.
\end{enumerate}
Moreover, any such representation is unique.
\end{coro}
\textbf{Proof.} The result is implied by Theorem \ref{finrepr} and by \cite{CB}, Theorem $1.2$ and \cite{CB}, Lemma $7.1$, $7.2$ once one observes that all finite subquivers of $Q$ are a disjoint union of Dynkin quivers.

\epf

 \subsection{The $SL(2,\C)$ case}\label{Khare}

As mentioned in Section \ref{infinitesimal}, we will now compare  Khare's result about representation theory of the deformed symplectic oscillator algebra of rank $1$ with the results of Section \ref{finitecase1} in the case $\G=SL(2,\C)$.  

We observe  that,  in this case, the subalgebra of invariant algebraic distributions supported at the identity can be identified with the algebra of  polynomials in the quadratic Casimir element $\Delta=\frac{1}{4}(EF+FE+\frac{H^2}{2})$ (where $E$, $F$, $H$ are the standard generators of $\mathfrak{sl}_2$), that coincides with the center of the enveloping algebra. Now if we let $x$, $y$ be a symplectic basis of the standard two dimensional complex symplectic vector space $V$, and we take $f=f(\Delta)$ to be a polynomial with no constant coefficient,  we can see that Khare's deformed symplectic oscillator algebra (cfr \cite{Kh}, section 9)
 $$
 H_f=\frac{TV\ltimes\mathcal{U}(\mathfrak{sl}_2)}{\langle [x,y]=1+f(\Delta)\rangle}
 $$ 
coincides with the infinitesimal Hecke algebra  $H_{c}(SL(2,\C))$  when we take $f=f_c$ to be an appropriate polynomial depending on $c$. 

Let   $V_C(i)$, $i\in\N_0$ be the standard cyclic module of $\mathfrak{sl}_2$  of highest weight $i$ (i.e. the irreducible finite dimensional representation of  dimension $i+1$).  Denote by  $b_i$ be the scalar by which the Casimir $\Delta$ acts on $V_C(i)$ ($b_i=i(i+2)/8$). 

Khare's classification of finite dimensional representations of $H_f$ can be summarized as follows (cfr \cite{Kh} Theorem $11$,  $\S$ $14$, and formula $(1)$, $\S$ $9$). 
\begin{enumerate}
 \item [I)] There exists a unique simple $H_f$-module of the form
\begin{enumerate}\item[($\ast$)]
$\displaystyle{\quad V(r,s):=\bigoplus_{i=s}^rV_C(i)}$
\end{enumerate}
for any $s\leq r$ in $N_0$ satisfying the two conditions: 
\begin{center}
\begin{enumerate}
\item[i)] $ \displaystyle{\quad \sum_{i=s}^{r}(i+1)(1+f(b_i))=0}$
\item[ii)] $\displaystyle{\quad\sum_{i=k}^{r}(i+1)(1+f(b_i))\neq 0 \qquad  s<k\leq r}$\, .
\end{enumerate}
\end{center}

\item[II)] Any finite dimensional irreducible $H_f$-module is isomorphic to one of the $V(r,s)$.  
\end{enumerate}

We observe now that all positive roots for the infinite quiver $A_{+\infty}$ are of the form $\alpha=\alpha_{[s,r]}=\sum_{i=s}^r \e_i$ for some $0\leq s \leq r$, where $\e_i$ are coordinate vectors (simple roots) as in Section \ref{roots} (cfr. \cite{K}, $\S$ $7.11$, where $\e_i=\alpha_i$ in Kac's notation). Thus according to Corollary \ref{CBHcont} in the case of $SL(2,\C)$ all possible simple, finite dimensional $\mathcal{H}_{c}(SL(2,\C))$ modules must have the form  $(*)$. Moreover, we observe that in this case $p(\beta)=0$ for any root, thus the conditions in Corollary  \ref{CBHcont} part $2)$ tell us that to have a representation of dimension vector $\alpha$ for  any decomposition of $\alpha=\beta^{(1)}+\cdots+\beta^{(n)}$ we must have $\lambda\cdot \beta^{(k)}=0$  for some $k$. Now, any decomposition looks like $\alpha_{[s,r]}=\alpha_{[s, s+t_1]}+\alpha_{[s+t_1+1,s+t_1+t_2]}+\cdots+\alpha_{[s+t_1+\cdots+t_{n}, r]}$. In particular, we can consider the decompositions 
 $\alpha_{[s, r]}=\alpha_{[s, m-1]}+\alpha_{[m, r]}$ for any $s<m\leq r$. Since $0=\lambda\cdot\alpha_{[s,r]}=\lambda \cdot \alpha_{[s, m-1]}+\lambda \cdot \alpha_{[m, r]}$ our condition implies that in particular $\lambda \cdot \alpha_{[m, r]}\neq 0$. On the other hand,  any nontrivial  decomposition of $\alpha_{[r,s]}$  contains a root $\alpha_{[m, r]}$ for $s<m\leq r$. Thus, the conditions in Corollary $\ref{CBHcont}$ can be rephrased as:
\begin{enumerate}
\item[a)] $\quad\alpha_{[s,r]}\cdot\lambda=0$
\item[b)] $\quad\alpha_{[m,r]}\cdot\lambda \neq 0\qquad s<m\leq r$.
\end{enumerate}

 We will now translate the conditions on  the dimension vector $\alpha$ in Corollary  \ref{CBHcont} part $(2)$ into Khare's conditions $i)$, $ii)$. In order to do this, we have to compare Khare's parameter  $f$ with our parameter $c$.  Let's  denote by $\chi_i$ the irreducible character corresponding to $V_C(i)$. Then , since the $\chi_i\, $s span the space of invariant functions,  we must have
\begin{equation}\label{sl2}
\lambda_i=\langle c,\chi_i\rangle=\langle 1+f(\Delta), \chi_i\rangle\, .
\end{equation}
  for any $i\in\N_0$. We recall now that constants in $\mathcal{U}(\mathfrak{sl}_2)$ correspond to multiples of the delta distribution $\delta_e$. Moreover for any $D\in \mathcal{U}(\mathfrak{sl}_2)$ one has  $D(\chi_i)(e)=\chi_i(D)=\tr_{V_C(i)}(D)$ , and, in particular $\tr_{V_C(i)}(\Delta^l)=\dim(V_C(i)) b_i^l=(i+1)b_i^l$. It is then easy to compute 
 \begin{eqnarray}
\langle 1+f(\Delta), \chi_i\rangle &=& (i+1)(1+f(b_i))\nonumber\,.
\end{eqnarray}
 and the equality (\ref{sl2}) becomes
\begin{equation}\label{equiparam}
(i+1)(1+f(b_i))=\lambda_i\, .
\end{equation}

Thus we can rewrite conditions $i)$ and $ii)$ as 
\begin{enumerate}
                                           \item[1)]$\displaystyle{\quad \sum_{i=s}^r\lambda_i=0} $
                                          \item[2)]$\quad \displaystyle{\sum_{i=k}^r\lambda_i\neq 0 \quad s<k\leq r}$
					  \end{enumerate}
which correspond to conditions $a)$, $b)$ above respectively.

\subsection{The higher rank case}

In  \cite{G}, W. L. Gan constructed reflection functors for higher rank. In the rank $1$ case reflection functors were constructed by Crawley-Boevey and Holland in \cite{CBH}.  Gan's reflection functors are defined for any loop-free vertex $i$ of any finite quiver $Q$ and, under some conditions on the parameter $\lambda, \nu$, they establish an equivalence between the categories of modules 
$$
F_i:\mathcal{A}_{n, \lambda, \nu}(Q)-\mathrm{mod}\rightarrow \mathcal{A}_{n, s_i\lambda, \nu}(Q)-\mathrm{mod}\, ,
$$
where $s_i\lambda$ denotes the action of the simple reflection $s_i$ at the vertex $i$ on the parameter $\lambda$. 
Thanks to this property the functors $F_i$ turned out to be a very powerful tool in the deformation-theoretic approach  to the study of finite dimensional representations of higher rank deformed preprojective algebras.

In this section we will briefly explain how Gan's definition of reflection functors can be pushed ahead, without any change, in the case of continuous deformed preprojective algebras of higher rank attached to an infinite affine Dynkin quiver and  finite dimensional representations. Gan's results on representations of the wreath product symplectic reflection algebra will then naturally extend to the continuous case. 

To ease notation, in  all that follows we will write $\Ann$ for  $\An$, and  $\hat{\Pi}_{\lambda}$ for $\Pc$, where $Q$ is an infinite affine quiver.

Let $i$ be any loop-free vertex of  $Q$. Since $\mathcal{A}_{n, \lambda, \nu}(Q)$ does not depend on the orientation of $Q$, we can suppose that $i$ is a sink (all  arrows at $i$ point toward $i$). Let $V$ be an $\mathcal{A}_{n, \lambda, \nu}(Q)$-module. Following Gan's construction we define  $F_i(V)$  as follows.

Let 
\begin{equation}\label{R}
R:=\left\{a\in Q| h(a)=i\right\}.
\end{equation}
Note that, for any infinite affine quiver, the set $R$ is  finite. 
If $\jj=(j_1,\dots,j_n)\in I^n$, where $I$ is the set of vertices of $Q$, let
$$
V_{\jj}:=|_{\jj}V, \qquad \mbox{\ and\ }\qquad \Delta(\jj):= \left\{ m\in[1, n]| j_m=i \right\}
$$
where $|_{\jj}$ is the element $e_{j_1}\otimes\cdots\otimes e_{j_n}$ as in Section \ref{highran}. For any subset $D\subset \Delta(\jj)$,  consider the finite set 
\begin{equation}\label{X}
\mathcal{X}(D):=\left\{\mbox{\ all maps\ } \xi: D\rightarrow R\right\}.
\end{equation}
For any $\xi\in\mathcal{X}(D)$ let 
$$
t(\jj, \xi):=(t_1,\dots,t_n)\in I^n,\qquad\mbox{where\ }\qquad t_m=\left\{\begin{array}{rl} j_m & \mbox{\ if \ }m\notin D\\ t(\xi(m)) & \mbox{\ if \ }m\in D \end{array}\right.
$$
Define 
$$
V(\jj, D):=\bigoplus_{\xi\in\mathcal{X}(D)}V_{t(\jj, \xi)}.
$$
so $V(\jj, \emptyset)=V_{\jj}$. For any $\xi$ there are  projection and  inclusion maps
$$
\pi_{\jj, \xi}:V(\jj, D)\rightarrow V_{t(\jj, \xi)}, \qquad \mu_{\jj, \xi}:V_{t(\jj, \xi)}\hookrightarrow V(\jj, D).
$$
Moreover, for any $p\in D$ there is a restriction map $\rho_p:\mathcal{X}(D)\rightarrow \mathcal{X}(D\setminus\{p\})$.
Thus for each $\xi\in\mathcal{X}(D)$ we can consider the two compositions
$$ \xymatrix{ V(\jj, D)\ar[rr]^{\pi_{\jj,\xi}} & & 
V_{t(\jj,\xi)} \ar[rr]^{\xi(p)_p|_{t(\jj,\xi)}} & & 
V_{t(\jj,\rho_p(\xi))} \ar@{^{(}->}[rr]^{\mu_{\jj,\rho_p(\xi)}}
& & V(\jj, D\setminus\{p\})}
$$
$$ \xymatrix{ V(\jj, D\setminus\{p\}) \ar[rr]^{\pi_{\jj,\rho_p(\xi)}}
& & V_{t(\jj,\rho_p(\xi))} 
       \ar[rr]^{\xi(p)^*_p\big|_{t(\jj,\rho_p(\xi))}} & & 
V_{t(\jj,\xi)} \ar@{^{(}->}[rr]^{\mu_{\jj,\xi}} & & V(\jj,D) }.
$$
Define now 
\begin{equation}\label{proj}
\pi_{\jj,p}(D):V(\jj, D)\rightarrow V(\jj, D\setminus\{p\}),\qquad \pi_{\jj,p}(D):=\sum_{\xi\in\mathcal{X}(D)}\mu_{\jj, \rho_p(\xi)}\xi(p)_p|_{t(\jj,\xi)}\pi_{\jj, \xi}
\end{equation}
\begin{equation}\label{incl}
\mu_{\jj,p}(D):V(\jj, D\setminus\{p\})\rightarrow V(\jj, D),\qquad \pi_{\jj,p}(D):=\sum_{\xi\in\mathcal{X}(D)}\mu_{\jj, \xi}\xi(p)^{\ast}_p|_{t(\jj,\rho_p(\xi))}\pi_{\jj, \rho_p(\xi)}.
\end{equation}
Let 
$$
V_{\jj}(D):=\left\{\begin{array}{rl} \cap_{p\in D}\Ker(\pi_{\jj, p}(D)) & \mbox{\ if \ }D\neq \emptyset\\ V_{\jj} & \mbox{\ if \ } D=\emptyset \end{array}\right.
$$
and let $V'_{\jj}:=V_{\jj}(\Delta(\jj))$. 
\begin{defi}\label{Fi}
We define $F_i(V):=V'=\bigoplus_{\jj\in I^n}V'_{\jj}$ as a vector space. 
\end{defi}
For any $l\in[1, n]$, $a\in \overline{Q}$, $\jj\in I^n$ with $j_l=t(a)$, we have to define maps $a'_{l}|_{\jj}:V'_{\jj}\rightarrow V'_{a_l(\jj)}$ , where  $a_l(\jj)$ is as defined in Section \ref{highran}. One has two cases:
\begin{itemize}
 \item[{\bf Case I.}]\  If $h(a),t(a)\neq i$ then $l\notin\Delta(\jj)=\Delta(a_l(\jj))$. For any $D\subset \Delta(\jj)$ we have a map
 $$
 a_l|_{\jj, D}:V(\jj, D)\rightarrow V(a_l(\jj), D),\qquad a_l|_{\jj, D}:=\sum_{\xi\in\mathcal{X}(D)}\mu_{a_l(\jj),\xi}a_l|_{t(\jj, \xi)}\pi_{\jj, \xi}.
 $$
We define
\begin{equation}\label{arrowone}
a'_l|_{\jj}:=a_l|_{\jj, \Delta(\jj)}.
\end{equation}

\item[{\bf Case II.}]\  If $t(a)=i$ then $l\in\Delta(\jj)$ and $\Delta(a_l(\jj))=\Delta(\jj)\backslash\left\{l\right\}$. For each $r\in R$ there is an injective map
$$
\tau_{r, l, D}:\mathcal{X}(D\backslash \left\{l\right\})\hookrightarrow \mathcal{X}(D):\eta\rightarrow \tau_{r, l, D}(\eta)
$$
where
$$
\tau_{r, l, D}(\eta)(m):=\left\{\begin{array}{rl} \eta(m) & \mbox{\ if \ }m\in D\backslash\left\{l\right\}\\ r & \mbox{\ if \ }m=l\end{array}\right.
$$
Since $t(\jj, \tau_{r, l, D}(D))=t(r^{\ast}_l(\jj), \eta)$, there is a projection map 
$$
\tau^{!}_{r, l, \jj, D}:V(\jj, D)\rightarrow V(r^{\ast}_l(\jj), D\backslash \left\{l\right\}),\qquad \tau^{!}_{r, l, \jj, D}:=\sum_{\eta\in\mathcal{X}(D\backslash \left\{l\right\})}\mu_{r^{\ast}_l(\jj), \eta}\pi_{\jj, \tau_{r, l, D}(\eta)},
$$
and an inclusion map
$$
\tau_{{r, l, \jj, D}_{!}}:V(r^{\ast}_l(\jj), D\backslash \left\{l\right\})\rightarrow V(\jj, D),\qquad \tau_{{r, l, \jj, D}_{!}}:=\sum_{\eta\in\mathcal{X}(D\backslash \left\{l\right\})}\mu_{\jj, \tau_{r, l, D}(\eta)}\pi_{r^{\ast}_l(\jj), \eta}.
$$
We define
\begin{equation}\label{arrowtwo}
a'_l|_{\jj}:=\tau^{!}_{a^{\ast}, l, \jj, \Delta_{\jj}}
\end{equation}
\item[{\bf Case III.}]\ If $h(a)=i$ then $l\notin\Delta(\jj)$ and $\Delta(a_l(\jj))=\Delta(\jj)\cup\left\{l\right\}$.
For any $D\subset \Delta(\jj)$ we have the inclusion map
$$
\tau_{{a, l, a_l(\jj), D\cup\left\{l\right\}}_{!}}:V(\jj, D)\rightarrow V(a_l(\jj), D\cup \left\{l\right\})
$$
as above. We have a map
$$
\theta_{a, l, \jj, D}:V(\jj, D)\rightarrow V(a_l(\jj), D\cup \left\{l\right\})
$$
defined by
$$
\theta_{a, l, \jj, D}:=\left(-\lambda_i+\mu_{a_l(\jj), l}\pi_{a_l(\jj), l}+\nu\sum_{m\in D}\s_{ml}|_{a_l(\jj)}\right)\tau_{{a, l, a_l(\jj), D\cup\left\{l\right\}}_{!}}.
$$
We define
\begin{equation}\label{arrowthree}
a'_l|_{\jj}:=\theta_{a, l, \jj, \Delta(\jj)}.
\end{equation}
\end{itemize}
 We define an action of $T_{\B}\E\rtimes\C[S_n]$ on $F_i(V)$, where $a_l|_{\jj}\in\E$ acts as $a'_l|_{\jj}$. We have the following proposition
\begin{prop}\label{Gan}{\rm(\cite{GG}, Proposition 2.7)}
 With the above action $F_i(V)$ is a $\mathcal{A}_{n, s_i\lambda, \nu}(Q)$-module.
\end{prop}
 
\epf

The so defined reflection functors for finite  dimensional representations of  the algebra $\Ann$  satisfy the same properties as in Gan (\cite{G}, Section $6.2$). In particular if $i$ is a loop-free vertex, let  $\Lambda_i$  be the set of all $(\lambda, \nu)\in B\times \C$ such that $\lambda_i\pm \nu\sum_{m=2}^{r}\s_{1m}$ is invertible in $\C[S_r]$ for  all $r\in [1,n]$.
\begin{thm}[\cite{G}, Theorem $5.1$]\label{catequi}
 If $(\lambda, \nu)\in \Lambda_i$, then the functor
 $$
 F_i: \Ann-mod\rightarrow \mathcal{A}_{n,r_i\lambda,\nu}-mod
 $$
is an equivalence of categories with quasi-inverse functor $F_i$.
\end{thm}
 \epf
\begin{rem}\label{lambda} Note that by \cite{G}, Proposition $5.12$ we have:
$$
\Lambda_i=\{(\lambda, \nu)\in B\times\C|\lambda_i\pm p\nu\neq 0 \mbox{\ for\ } p=0,\dots, n-1\}
$$ 
\end{rem}

 From now on we will consider $\nu, \lambda$ as formal parameters. More specifically, let $U$ be a finite dimensional complex vector space  and consider $\nu\in\K:=\C[[U]]$, $\lambda\in B=\Pi_{i\in I}\K$. 
Let  $\mathfrak{m}$ be the unique maximal ideal of $\K$.  For any $\K$-module write $\overline{V}:=V/\mathfrak{m}V$.
A $\Ann$-module $V$ is a flat deformation of a $\overline{\Ann}$-module $V_0$ if $V\cong V_0[[U]]$ as $\K$-modules and $\overline{V}\cong V_0$. 
Consider  the decreasing filtration $\Ann\supset \mathfrak{m}\Ann\supset\mathfrak{m}^2\Ann\supset\dots$\ . We have that  $\mathrm{Gr}_{\mathfrak{m}}\Ann\cong\overline{\Ann}[[U]]$ as algebras over $\K$.  Then Lemma $5.13$ and Proposition $5.14$ of \cite{G} extend to the continuous case with analogous proofs. Let $\tilde{B}=\prod_{i\in I}\C[U]$.

\begin{lem}\label{defor1}
 Assume $\lambda\in\tilde{B}$ and $\nu\in \C[U]$. Then $\mathrm{Gr}_{\mathfrak{m}}\Ann\cong\overline{\Ann}[[U]]$ as algebras over $\C$. 
\end{lem}

\epf

\begin{prop}\label{defor2}
 Assume $\lambda$ and $\nu$ are as in Lemma $\ref{defor1}$. Let $i\in I$ and suppose $\lambda_i\pm\nu\sum_{m=2}^r$ is invertible in $\mathbb{K}[S_r]$, for any $r\in [1, n]$. If a $\Ann$-module $V$ is a flat deformation of a $\overline{ \Ann}$-module $V_0$, then  $F_i(V)$ is a flat deformation of $F_i(V_0)$. 
\end{prop}

\epf

Let now $\vec{n}=(n_1,\dots,n_r)$ be a partition of $n$. $X=X_1\otimes\cdots\otimes X_r$ be a simple module for the group $S_{\vec{n}}:=S_{n_1}\times\cdots\times S_{n_r}\subset S_n$ and let $\{ i_1,\dots,i_r\}$ be $r$ distinct vertices of $Q$. For any $i\in I$, let $\Nn_i$ be the complex vector space with dimension vector $\e_i$ and let $\Nn=\Nn_{i_1}^{\otimes n_1}\otimes\cdots\otimes \Nn_{i_r}^{\otimes n_r}$. Then $X\otimes \Nn$ is a simple module for $\Bb\ltimes \C[S_{\vec{n}}]$. One can then form the induced module $X\otimes \Nn\uparrow:=\mathrm{Ind}_{\Bb\ltimes \C[S_{\vec{n}}]}^{\Bb\ltimes\C[S_n]}(X\otimes \Nn)$ over $\Bb\ltimes\C[S_n]$. Moreover, it is known that any simple finite dimensional  $\Bb\ltimes\C[S_n]$-module has this form (this is true by \cite{Mac}, paragraph after (A5), when $\G$ is finite, and remains true for $\G$ reductive when we consider only finite dimensional representations).

The next Lemma is the analog of Lemma $6.1$ in \cite{G}  and it can be proved in the same way. 
\begin{lem}
 Let the $\Ann$-module $V$ be a flat deformation of  the $\Ann$-module $\overline{V}$. If $\overline{V}$ is simple as a $\Bb\ltimes \C[S_{n}]$-module, then all elements of $\E$ must act by $0$ on $V$.
\end{lem}

\epf
 
The following Theorem can be proved as  Theorem $6.2$ in \cite{G} and is equivalent to Theorem $6.5$ in \cite{EGG}. 
\begin{thm}\label{exten1}
 Assume $\nu\neq 0$. The $\B\rtimes \K[S_n]$-module $(X\otimes\Nn\uparrow)[[U]]$ extends to a $\Ann$-module if and only if the following conditions are satisfied:
 \begin{itemize}
 \item [(i)] For all $l\in[1,\dots,r]$, the simple module $X_l$ of $S_{n_l}$ has rectangular Young diagram of size $a_l\times b_l$;
\item [(ii)] no two vertices in the collection $\{i_1,\dots,i_r\}$ are adjacent in $Q$, i.e. $\langle \e_{i_j}, \e_{i_k}\rangle=0$ for any $j\neq k$;
 \item[(iii)] for all $l\in [1,\dots, r]$, one has $\lambda_{i_l}=\nu(a_l-b_l)$. 
 \end{itemize}
Where we agree that condition $(ii)$ is empty if $r=1$, that is to say in the case of the trivial partition $\vec{n}=(n)$.
\end{thm}

\epf

 Let now $\lambda_0\in\overline{B}$ ($B=\hat{\K}^I$). We will write $\Pco$ for $\overline{\Pco}$ and $\Ano$ for $\overline{\Ano}$.
 
  Let $Y_1,\dots,Y_r$ be a collection of pairwise non-isomorphic representations of $\Pco$ with dimension vectors $\beta_{(1)},\dots,\beta_{(r)}$ respectively. Let $Y:=Y_1^{\otimes n_1}\otimes\cdots\otimes Y_r^{\otimes n_r}$. Then $X\otimes Y$ is an irreducible representation of $\Pco^{\otimes n}\rtimes\C[S_{\vec{n}}]$ and, as before, we can consider the induced $\Ano$-module $X\otimes Y\uparrow:=\ind_{\Pco^{\ottimes n}\rtimes\C[S_{\vec{n}}]}^{\Pco^{\ottimes n}\rtimes\C[S_n]}X\otimes Y$. It is known (as before by \cite{Mac}) that any finite dimensional simple $\Ano$-module is of this form.

Now  for any $\hat{\Pi}_{\lambda}$-module $M$ we define
 
 $$
 \Su(M):=\{ i\in I| e_iM\neq0\}.
 $$   

Similarly, for any $\Ann$-module $V$ we define 
 
 $$
 \Su(V):=\{ i\in I| \exists  \,  \jj=(j_1,\dots, j_n) \mbox{\ such that\ } V_{\jj}\neq 0 \mbox{\ and\  }  i\in\left\{ j_1,\dots, j_n\right\}\}.
 $$   

\begin{rem}\label{support} 
We observe the following facts:
\begin{itemize} 
\item we have 
$$
\Su(X\otimes\Nn \uparrow)=\Su(X\otimes \Nn)=\{i_1,\dots,i_r\}^{n}
$$ 
$$
\Su(X\otimes Y\uparrow)=\Su(X\otimes Y)=\left(\cup_{i=1}^{r}\Su(Y_i)\right)^{n}\, ;
$$ 
\item let $J\subset I$ be the minimal subset of vertices such that $J$ corresponds to a connected subquiver $Q_J$ of $Q$ and $\cup_{i=1}^{r}\Su(Y_i)\subset J$.  From the proof of Lemma \ref{dominance} we know it exists an element $w\in W_J\subset W$  such that $w\lambda_0=\lambda^{+}$, with $\lambda^{+}$ $J$-dominant and $w$ of minimal length as an element of $W(Q_J)$. If   $w=s_{j_m}\cdots s_{j_1}$ is a reduced expression for $w$, where $m$ is the length of $w$ in $W(Q_J)$, then by minimality we have $r_{j_{k-1}}\cdots r_{j_1}\lambda_0\cdot \e_{j_{k}}\neq 0$ for all $k\in [1,\cdots,m]$.  We will write $F_w=F_{j_m}\cdots F_{j_1}$;
 \item from the definition of the reflection functor $F_i$ it is clear that if $i\in\Su(V)$ then $\Su(F_i(V))\subset\Su(V)$. In particular we have 
$$
\Su \left(F_w(X\otimes Y\uparrow)\right) \subset \Su(X\otimes Y\uparrow) ;
$$
\item from \cite{CBH} Lemma $7.1$, we can deduce that if $M$ is an irreducible representation of $\hat{\Pi}_{\lambda^+}$ and $\Su(M)\subset J$ then $\Su(M)=\{j\}$ for some $j\in J$ with $ \lambda^+\cdot \e_j=0$. In particular, if we consider  the set $\Sigma_{\lambda_0, J}$ of dimension vectors of the representations of $\Pco$ with support in $J$, we have $w\Sigma_{\lambda_0, J}=\Sigma_{\lambda^+, J}=\{\e_j | j\in J,\quad(\lambda^+)_j=0\}$.
Thus, if $\beta_{(l)}$ is the dimension vector of $Y_l$, for any $l\in[1,\dots r]$ we have that  $w\beta_{(l)}=\e_{i_l}$ for some  $i_l\in J$ with $\lambda_{i_l}=0$.  As a consequence we must have  $\Su~(~F_w~(~X~\otimes~ ~Y~\uparrow~)~)~=~\{i_1,\dots, i_r\}^n\subset (\Sigma_{\lambda^{+}, J})^n$.
\end{itemize}
\end{rem}

 The following theorem is the analog of  Theorem $6.3$ in \cite{G}, and of Theorem $1.3$ in \cite{M}.
 
 \begin{thm}\label{mainhigh}
Let $\lambda\in B$. Assume $\lambda_i\in U$ for all $i$ and $0\neq\nu\in U$. The $\Ano$-module $X\otimes Y\uparrow$ has a flat deformation to a $\Anon$-module if and only if the following conditions are satisfied:
\begin{enumerate}
 \item [(i)] for all $l\in[1,\dots,r]$, the simple module $X_l$ of $S_{n_l}$ has rectangular Young diagram, of size $a_l\times b_l$;
 \item[(ii)] $\langle \beta^l, \beta^m\rangle=0$ for all $l\neq m$;
 \item[(iii)] for all $l\in [1,\dots,r]$, one has $\lambda\cdot \beta_{(l)}=(a_l-b_l)\nu$, where $\beta_{(l)}$ is the dimension vector of $Y_l$\,;
 \end{enumerate}
 where we agree that condition $(ii)$ is empty  in the case of the trivial partition $\vec{n}=(n)$.
 \end{thm}

 \pf\  The proof goes exactly as  in \cite{G} Theorem $6.3$ when one observes that, by Remark \ref{support} above, one has $F_w(X\otimes Y\uparrow)= X\otimes \Nn\uparrow$, where $i_l=w\beta_{(l)}$.
 
\epf
\begin{rem} It can be deduced from \cite{CBH}, $\S$ $7$,  that the condition $\langle\beta_{(i_l)}, \beta_{(i_m)}\rangle=0$ for the two distinct roots $\beta_{(l)}$, $\beta_{(m)}$ is equivalent to the condition: $\Ext^1_{\hat{\Pi}_{\lambda_0}}(Y_l, Y_m)=0$ for the two irreducible non-isomorphic representations $Y_l$, $Y_m$ of $\hat{\Pi}_{\lambda_0}$ (this is true when $\lambda_0$ is $J$-dominant  by \cite{CBH},  Lemma $7.1$, and can be proved  for any weight using the reflection functors  and our Lemma $\ref{dominance}$). The last one is the form condition $(ii)$ is stated in  both \cite{M}  and \cite{G}.  Note that when $Y_l\cong Y_m$ the two conditions are not equivalent anymore. Indeed, by Corollary $7.6$ in \cite{CBH} (see also \cite{EM}, Proposition $4.6$), we have  $\Ext^1_{\hat{\Pi}_{\lambda_0}}(Y, Y)=0$ for any irreducible finite dimensional representation $Y$ of $\hat{\Pi}_{\lambda_0}$, while the Ringel form attached to any infinite affine Dynkin quiver is positive definite  on the set of roots (any principal minor of the corresponding cartan matrix is positive definite and any root is a linear combination of a finite number of simple roots). For this reason we specified condition $(ii)$ is empty in the case of the trivial partition.
\end{rem}

Let now $\lambda$, $\nu$ be regular functions on the vector space $U$ such that condition $(iii)$ of Theorem \ref{mainhigh} holds. Suppose there exists a point $o\in U$ such that $\lambda$ specializes  to $\lambda_0$ and $\nu=0$ at $p$. Let $j_m,\dots, j_1$ be as in Remark \ref{support}. According to Gan's notation  (\cite{G}, $\S$, $6.3$) let $U'$ be the Zariski open set in $U$ defined by $(r_{j_k}\dots r_{j_i}\pm p\nu\neq 0$ for all $k=[1,\dots ,m]$ and $p=0,\dots, n-1$ ( see Remark \ref{lambda}). Clearly $o\in U'$.  Let $\C[U']$ be the algebra of regular functions on $U'$ and, for any $u\in U'$, let $\mathfrak{m}_{u}$ be the maximal ideal of  functions vanishing at $u$. For any $\C[U']$-module $V$ write $V^u=V/\mathfrak{m}_u V$. We have the following theorem the proof of which is as in \cite{G}, Theorem $6.4$.
\begin{thm}\label{specialization}
There exists a $\Ann$-module $V$ such that:
\begin{enumerate}
\item[(i)] $V^{o}=X\otimes Y\uparrow$ as a $\Ann$-module and $V$ is flat over $U'$; 
\item[(ii)] for any point $u\in U'$, $V^{u}$ is a finite dimensional simple $\Ann^{u}$-module, isomorphic to $X\otimes Y\uparrow$ as a $\mathcal{B}^u\ltimes \C[S_n]$-module.
\end{enumerate}
\end{thm}

\epf

\textbf{Acknowledgments}. I am very grateful to Pavel Etingof for many comments and useful discussions and to Wee Liang Gan for  patiently explaining to me  the results of \cite{GG}, \cite{G}.

\end{document}